\documentclass[11pt]{article}
\usepackage[utf8]{inputenc}
\usepackage{color}
\usepackage{xcolor}
\usepackage{ulem}

\usepackage{amsmath,amsfonts,euscript,amssymb,amsthm,a4,times}

\allowdisplaybreaks

\def\XXint#1#2#3{{\setbox0=\hbox{$#1{#2#3}{\int}$}
\vcenter{\hbox{$#2#3$}}\kern-.5\wd0}}

\newtheorem{theorem}{Theorem}[section]

\newtheorem{lemma}[theorem]{Lemma}

\newtheorem{proposition}[theorem]{Proposition}

\theoremstyle{definition}
\newtheorem{definition}[theorem]{Definition}
\newtheorem{remark}[theorem]{Remark}

\def\tcr{\textcolor{red}} 
\def\tcb{\textcolor{blue}} 

\makeatletter \catcode`@=11
\newbox\tr@tto
\setbox\tr@tto=\hbox{{\count0=0\dimen0=-,9pt\dimen1=1,1pt\loop\ifnum
    \count0<11 \advance \count0 by1 \vrule width.51pt height\dimen1
    depth\dimen0\kern-0.17pt\advance\dimen0 by-0.05pt\advance\dimen1
    by0.1pt\repeat \loop\ifnum\count0<21\advance \count0 by1 \vrule
    width.6pt height\dimen1 depth\dimen0\kern-0.2pt \advance\dimen0
    by-0.1pt\advance\dimen1 by 0.05pt\repeat}}
\def\medint{\displaystyle\copy\tr@tto\kern-10.4pt\int}
\numberwithin{equation}{section}

\def\R{{\mathbb R}}

\def\proofof#1{\begin{proof}[Proof of #1]}

\def\dx{{\mathrm d}x}

\def\loc{\rm loc}

\def\R{{\mathbb R}}

\def\proofof#1{\begin{proof}[Proof of #1]}

\def\dx{{\mathrm d}x}

\def\loc{\rm loc}

\catcode`@=11
\newbox\tr@tto
\setbox\tr@tto=\hbox{{\count0=0\dimen0=-,9pt\dimen1=1,1pt\loop\ifnum\count0<11 \advance \count0 by1 \vrule width.51pt height\dimen1
     depth\dimen0\kern-0.17pt\advance\dimen0 by-0.05pt\advance\dimen1
     by0.1pt\repeat \loop\ifnum\count0<21\advance \count0 by1 \vrule
     width.6pt height\dimen1 depth\dimen0\kern-0.2pt \advance\dimen0
     by-0.1pt\advance\dimen1 by 0.05pt\repeat}}
\def\medint{\displaystyle\copy\tr@tto\kern-10.4pt\int}
\catcode`@=12

\numberwithin{equation}{section}

\date{\today}

\begin{document}


\title{Non autonomous problems with {nearly} linear growth and Orlicz-Sobolev coefficients}
\author{Michela Eleuteri $\&$ Antonia Passarelli di Napoli}

\maketitle

\begin{abstract}
We establish the local Lipschitz regularity of the local minimizers of non autonomous integral {functionals} of the form
\begin{equation*}
\label{functional}
{\mathcal{F}({u, \Omega}) =} \int_\Omega F(x, {Du})\, \dx,
\end{equation*}
where $\Omega$ is a bounded open set of $\mathbb{R}^n$, $n \ge 2$ and the energy density $F(x,\xi)$ satisfies {nearly} linear growth conditions with respect to the gradient variable.
The main novelty of this work lies in the identification of a {suitable quantitative balance between the logarithmic degeneration of the ellipticity and an} Orlicz-Sobolev regularity of the partial map $x\to {F(x,\xi)}$  that is sufficient to obtain the Lipschitz regularity of the local minimizers. 

\end{abstract}



\bigskip
\section{Introduction}

{{\it Who of us would not be glad to lift the veil behind which the future lies hidden; to
cast a glance at the next advances of our science and at the secrets of its development
during future centuries?}}

\vspace{4mm}

{
With these words Hilbert opened his lecture at the International Congress of Mathematicians at Paris in 1900 \cite{H00}. Among the very famous list of open problems he presented, the 19th and 20th were devoted to the regularity theory.}
In the first half of the last century it was clear that every sufficiently smooth solution to a linear or nonlinear analytic elliptic system is an analytic function (for more details see \cite[Chapter 2, Section II]{Giaquinta}); however, in general, since it is possible to find extremals or weak solutions in some Sobolev spaces, a gap in the regularity scale appears.
 The pioneering work by De Giorgi in 1957 \cite{DeGio} filled this gap,  opening the way to the regularity of solutions to elliptic equations with bounded measurable coefficients and for minima of the corresponding functionals in the Calculus of Variations (see also \cite{Nash}, \cite{Moser}).  De Giorgi's theorem was later generalized in several directions (see for instance the book by Ladyzhenskaya and Ural'tseva \cite{LU}), in particular it was soon recognized that the linearity of the equation played actually no role in the proof.
 
{It took almost  twenty years to start exploiting the full impact of De Giorgi's
techniques on the regularity of minima until Giaquinta and Giusti in 1982 \cite{GiaquintaGiusti82} applied De Giorgi's method in a direct way, i.e. without passing
through the Euler-Lagrange equation. Subsequently a wide theory has been developed for uniformly elliptic variational integrals with {\it standard  growth}, namely for functionals of type 
\begin{equation}\label{functional}
    \mathcal{F}(u,\Omega)= \int_\Omega F(x,Du)\, \dx
 \end{equation} 
where 
\[
|\xi|^p\lesssim F(x,\xi)\lesssim |\xi|^p \qquad p > 1.
\]}
In this classical framework, minimizers enjoy many  regularity properties, such as local boundedness, higher integrability of the gradient, H\"older continuity of the first derivatives and  their higher differentiability. We refer to the classical monographs by Giaquinta \cite{Giaquinta} and Giusti \cite{Giusti} for an exhaustive treatment  of these developments.

Over the last three decades, increasing attention has been devoted to variational problems whose ellipticity is no longer uniformly controlled {(see for instance \cite{DFM23} together with the references therein).} Besides posing challenging mathematical questions, these models naturally arise in several applications, including nonlinear elasticity, plasticity, composite materials, image restoration, and the theory of electrorheological fluids. In these situations the ellipticity of the integrand may be lost as the gradient becomes large, producing phenomena that do not have a counterpart in the uniformly elliptic setting. {We refer to \cite{DF26} for a unified and comprehensive survey concerning the recent advances of the nonuniformly elliptic setting, with a particular reference to the basic model examples capturing some of the main aspects of the uniform/nonuniform ellipticity notions.}

{One of the prominent examples of rate of non uniformity for equations and functionals is the polynomial
one, i.e. when the pointwise ellipticity ratio behaves (at infinity) as a positive
power of the gradient variable. Since the foundational works of regularity {(see for instance \cite{Ser69,LU70,Iva72,Sim76,Tru67,Tru81,UU84})}, many significant progresses have been made, in particular after the innovative approach introduced by Marcellini}
 {who {was the first to start a systematic study of functionals with $(p,q)-$growth conditions in a series
of seminal papers \cite{M89,M91,M93}.} The prototype consists of functionals whose energy density satisfies
\[
|\xi|^p\lesssim F(x,\xi)\lesssim |\xi|^q,
\qquad
1<p\le q
\]
and the key condition 
is controlling the lowest and highest eigenvalues of $F_{\xi \xi}$
between two different powers of the gradient variable in a way to extend the $p-$Laplace
type behaviour, i.e.,
\[
\left \{
\begin{array}{lll}
\!\!\!\!\!\! & |{\xi}|^{p-2} \mathbb{I}_{n \times n} \lesssim {F_{\xi \xi}(x, \xi)} \lesssim |{\xi}|^{q-2} \mathbb{I}_{n \times n} \qquad \textnormal{in the sense of matrices}\\[3mm]
\!\!\!\!\!\! & \textnormal{for some $1 < p \le q$ and all $|{\xi}| \ge 1,$ $x \in \Omega$}
\end{array}
\right.
\]
so that the ellipticity ratio of the corresponding Euler Lagrange equation} {can be controlled in terms of the difference $q-p$}
\begin{equation}
\label{ellipticityratio}
{\mathcal R(x, {\xi})\approx |{\xi}|^{q-p}}, \qquad {|{\xi}| \ge \, 1,}
\end{equation}
thus becoming unbounded whenever $q>p$. Marcellini's pioneering papers revealed that the admissible gap $q-p$ is strictly connected with the regularity of minimizers, opening a research direction that has {deeply} influenced the modern regularity theory.

Since then, remarkable progress has been achieved. Sharp or nearly optimal regularity results have been established under increasingly weaker assumptions on the ellipticity ratio, the dimension, and the regularity of the coefficients. We refer, among many others, to the contributions \cite{BS1, BS2, CKP, CKPPisa, CKPComm, CGG, Choe92,  CGGP0, CGGP, CGHPdN20, CM15, CM18, CM23, CM24, CM25,  DFM, DFM2, FDFL,   EMM16, EMM2, EMM3, EMMP, EPT, ELM1, ELM2, EM,  M2021, Mar21, M2022, M06, Mingione-Radulescu}
and the references therein.

 The regularity theory developed over the last decade {therefore} strongly indicates that the precise  structure of the energy is often less important than the asymptotic behaviour of the ellipticity ratio itself, {which reveals to be} the fundamental quantity governing the regularity of minimizers.

A second fundamental aspect concerns the dependence of the energy density on the spatial variable. In the autonomous case,
$
F(x,\xi)\equiv F(\xi),
$
the regularity of minimizers is entirely determined by the {structure} and the ellipticity of the integrand. Once the dependence upon $x$ is introduced, however, the regularity of the coefficient map
$
x\mapsto F(x,\xi)
$
becomes an  equally crucial ingredient.

Already in the uniformly elliptic setting it is well known that quantitative regularity estimates cannot be expected without suitable assumptions on the coefficients. During the last few years, considerable effort has therefore been devoted to identifying the weakest possible regularity assumptions ensuring local Lipschitz continuity and higher differentiability of minimizers. This line of research has progressively moved from H\"older continuous coefficients to Sobolev and Orlicz--Sobolev regularity assumptions, leading to several optimal and borderline results; see, for instance,
\cite{CGHPdN20, CMMPdN, DFM,  DFM2,  EMM16, EMM2, EPdN23, EPdN24}
and the references therein.

Despite these significant advances, the simultaneous interaction between the non-uniform ellipticity and low regularity of the coefficients is still far from being completely understood. Indeed, these two phenomena act in opposite directions. On one hand, the loss of ellipticity tends to deteriorate the regularity of minimizers; on the other hand, additional regularity of the coefficients can compensate, at least partially, for this degeneration. 

A systematic investigation of this interplay has recently been initiated in
\cite{EMMPdN26, EPdN}, where quantitative conditions relating the growth of the ellipticity ratio to the Sobolev regularity of the coefficients have been identified. The present paper continues this line of research in a substantially different setting, namely that of variational integrals exhibiting {{nearly} linear growth or {almost linear growth in the terminology of \cite{GIS}}. As we shall explain below, this regime falls outside the currently available theory and requires new ideas, owing to the borderline character of both the growth conditions and the natural energy space.
 More precisely, we investigate functionals of the form
{\eqref{functional}}
whose energy density exhibits {almost linear growth}. More precisely, we assume that
\[
\lim_{|\xi|\to\infty}\frac{F(x,\xi)}{|\xi|}=+\infty,
\qquad
\lim_{|\xi|\to\infty}\frac{F(x,\xi)}{|\xi|^p}=0
\quad\text{for every }p>1.
\]
Thus, the growth of the integrand is only {slightly} stronger than linear and remains below every polynomial growth of order strictly larger than one.

Although this regime may at first sight appear as a limiting case of the theory of $(p,q)$-growth, its analytical features are largely different. Indeed, the natural energy space is no longer a Sobolev space $W^{1,p}$, but rather an Orlicz--Sobolev space generated by an $N$-function with {almost} linear growth. Consequently, several fundamental tools that are usually employed in the polynomial setting, including scaling arguments and Sobolev embeddings,  are no longer available in their standard form.

The model case motivating our analysis is provided by integrands {of the form}
\[
F(x,\xi)\approx |\xi|\log(1+|\xi|),
\]
whose natural energy space is \(
W^{1,L\log L}_{\mathrm{loc}}(\Omega).
\)
{(See Subsection 2.1 for the definition and properties).}
Assumptions of this type still generate non-uniformly elliptic Euler equations. Indeed, under the structural assumptions introduced below, the eigenvalues of the Hessian matrix
{$F_{\xi \xi}(x,\xi)$}
satisfy
\begin{equation}
\label{ratioeigen}
\frac{\Lambda_{\max}(\xi)}
     {\Lambda_{\min}(\xi)}
\approx
\log(1+|\xi|),
\end{equation}
so that the ellipticity ratio blows up as $|\xi|\to \infty$, although only at a logarithmic rate. This behaviour should be compared with the polynomial estimate {\eqref{ellipticityratio}}
arising in the $(p,q)$-growth theory. While the degeneration is considerably milder, it occurs precisely when the  coercivity is also weakest. Consequently, neither the techniques developed for uniformly elliptic integrals nor those for polynomially non-uniformly elliptic problems apply directly.

The autonomous counterpart of this class of functionals has been extensively investigated during the last decades. Starting from the pioneering work of {Greco, Iwaniec and Sbordone} \cite{GIS} , higher integrability and regularity results have been established for a broad class of almost linearly growing energies. Subsequently, partial and full $C^{1,\alpha}$ regularity results have been obtained under increasingly general assumptions; see, for instance,
\cite{BF01, BDPdN, DDP, DDH, DM23, DP24, EM, FM00, FO98, FS98, MS99}
and the references therein.

The non-autonomous case is considerably less developed. Existing regularity results  assume that the dependence upon the spatial variable is sufficiently regular, typically H\"older continuous. Under a  Orlicz-Sobolev regularity assumptions stronger than the borderline space $W^{1,n}$, Giova in \cite{Giova} established an higher differentiability result but,
to the best of our knowledge, no Lipschitz regularity results are  available under such a borderline Orlicz--Sobolev assumptions on the coefficients.

The purpose of this paper is  to fill this gap. Our {aim} is to identify a {suitable} quantitative balance between the logarithmic degeneration of the ellipticity and the minimal regularity of the coefficients ensuring local Lipschitz continuity of minimizers. In this respect, our work should be regarded as the natural continuation, in the almost linear setting, of the program initiated in
\cite{EMMPdN26,EPdN} for more general classes of non-uniformly elliptic functionals.
\\
To state our assumptions, let
$
F:\Omega\times\mathbb R^n\rightarrow[0,+\infty)
$
be a Carath\'eodory function such that, for almost every
$x\in\Omega$, the map
$\xi\mapsto F(x,\xi)$
belongs to {$C^2(\mathbb R^n)$ and that $F(x,0)=0$.}
Moreover, we assume that the dependence upon the gradient is radial, {that is}
$
F(x,\xi)=\widetilde F(x,|\xi|),
$
for a suitable function
$\widetilde F:\Omega\times[0,\infty)\rightarrow[0,\infty)$ and that {$F_{\xi x}(x,\xi)$} exists and it is  a Carath\'eodory function.

The structural assumptions imposed throughout the paper are the following.
First, we require
{that the energy density has almost linear growth,
more precisely
\[
c|\xi|\log(1+|\xi|)
\le
F(x,\xi)
\le
C(1+|\xi|\log(1+|\xi|)) \tag{F1},
\]
for $0 < c < C$ and moreover we also ask} that the Hessian matrix with respect to the gradient variable satisfies the ellipticity bounds
\[
\lambda\frac{|\mu|^2}{1+|\xi|}
\le
\langle {F_{\xi\xi}}(x,\xi)\mu,\mu\rangle
\le
\Lambda
\frac{\log(1+|\xi|)}{|\xi|}
|\mu|^2,
\tag{F2}
\]
{for every $\xi, \mu \in \mathbb{R}^n$ and a.e. $x$ in $\Omega,$}
where {$0 < \lambda < \Lambda$ and where} we supposed that $\frac{\log(1+t)}{t}_{|t=0}=1$. 

The dependence upon the spatial variable is measured through the following inequality
\[
|{F_{\xi x}}(x,\xi)|
\le
h(x)\log(1+|\xi|),
\tag{F3}
\]
where $h$ is a non-negative measurable function, {satisfying assumption {(F4)} below.}

{These assumptions show that} the corresponding Euler equation {associated to our functional \eqref{functional}} is genuinely non-uniformly elliptic.
The coexistence of a borderline coercivity with a logarithmic degeneration represents the main analytical feature of the class of functionals considered here.

The dependence of the coefficients upon the spatial variable is encoded by the function $h$ appearing in (F3). In order to compensate for the loss of ellipticity, we assume the borderline Orlicz integrability condition
\[
\int_{B_R}
h^n
\log^{\alpha n}(1+h)\, \dx
<
+\infty,
\qquad
B_R\Subset\Omega,
\tag{F4}
\]
where
\[
\alpha>
\begin{cases}
4, & n\ge3,\\[2mm]
\displaystyle {\frac92}, & n=2.
\end{cases}
\]

Assumption {(F4)} means precisely that
\[
h\in
L^n\log^{\alpha n}L_{\mathrm{loc}}(\Omega),
\]
{(see Section \ref{Orlicz} for more details)}} which is only slightly stronger than the critical space
$L^n$.
This should be compared with the recent regularity theory for non-uniformly elliptic problems with polynomial growth, where the borderline threshold is also given by a suitable logarithmic improvement of the critical Sobolev regularity; see
\cite{CGHPdN20,EPdN23,EMMPdN26}.
Therefore, condition {(F4)} appears as the natural counterpart of those assumptions in the present almost linear framework.

We emphasize that the class of coefficients covered by {(F4)} is substantially larger than the H\"older continuous one considered in the previous literature on almost linearly growing energies,  while remaining sufficiently regular to compensate for the logarithmic degeneration of the ellipticity.

Before stating our main result, we give the notion of local minimizer.

\begin{definition}
We say
that $u\in W_{loc}^{1,1}(\Omega )$ is a \textit{local minimizer} of the
integral functional $\mathcal{F}$ in \eqref{functional} if $F(x, Du) \in
L_{loc}^{1}(\Omega )$ and 
\begin{equation*}
\int_{\Omega ^{\prime }}F(x, Du)\, \dx\leq \int_{\Omega ^{\prime }}F(x, Du+D\varphi
)\, \dx
\end{equation*}%
for every open set $\Omega ^{\prime }$, $\overline{\Omega ^{\prime }}\subset
\Omega $ and for every $\varphi \in W_{0}^{1,1}(\Omega ^{\prime })$. 
\end{definition}

\color{black}

We now state the main result of the paper.

\begin{theorem}
Let {$u\in W^{1,1}_{\mathrm{loc}}(\Omega)
$} 
be a local minimizer of the functional at \eqref{functional}
under assumptions {{\rm (F1)--(F4)}}.
Then
\[
u\in W^{1,\infty}_{\mathrm{loc}}(\Omega),
\]
and the following estimate
\begin{eqnarray}\label{estfin}
\|Du\|_{L^\infty(B_\rho)}
&\le&
C\int_{B_{  R}}\Big(1+|Du|\log(1+|Du|)\Big)
\, \dx\cr\cr 
&&\qquad\cdot \exp\Bigg(C\int_{B_{ R}}\Big(1+|Du|\log(1+|Du|)\Big)\, \dx\Bigg)\cr\cr
    &&+\frac{C}{(  R- \rho )^{{\theta}}}\int_{B_{ R}}\Big(1+|Du|\log(1+|Du|)\Big)\, \dx\cr\cr
&&\qquad\cdot\log^{{\theta}}\left(1+\frac{C}{(  R- \rho )^{{\theta}}}\int_{B_{  R}}\Big(1+|Du|\log(1+|Du|)\Big)\, \dx\right) \cr\cr
 &&+C\left(\frac{1}{(  R- \rho )^{{\theta}}}+1\right)\int_{B_{ R}}\Big(1+|Du|\log(1+|Du|)\Big)\, \dx
\end{eqnarray}
holds  for every pair of concentric balls
$
B_\rho\Subset B_R\Subset\Omega,
$ with a constant
$C=C(n,\lambda,\Lambda, h)$, {where \[
\theta =
\begin{cases}
n, & n\ge3,\\[2mm]
2q, & n=2,
\end{cases}
\]
for any $q < 1.$}
\end{theorem}

The estimate above provides a local Lipschitz bound depending only on the natural energy of the minimizer and on the structural constants of the problem. In particular, it shows that the logarithmic improvement of the regularity assumption on the coefficients is sufficient to compensate for the logarithmic degeneration of the ellipticity.

{As usual in the regularity theory of variational integrals,} the proof combines a suitable approximation argument with a new {a priori} estimate for smooth minimizers.

Actually, the original functional is first approximated by a family of uniformly regular functionals admitting unique smooth minimizers. Uniform estimates are then established for the approximating solutions and finally transferred to the original minimizer through a compactness argument.
\\
The crucial point consists in deriving the uniform {a priori} estimate for smooth minimizers under the borderline assumptions {{\rm (F1)--(F4)}}.
In this step, the main difficulty relies in the    estimate of the terms generated by the dependence of the coefficients on the spatial variable without losing control of the logarithmic degeneracy of the ellipticity.

The key idea here is to exploit the additional logarithmic summability of the coefficient function $h$ in assumption {{\rm (F4)}} in a quantitatively optimal way. More precisely, the extra integrability encoded by the Orlicz space
$
L^n\log^{\alpha n}L
$
is shown to exactly balance the logarithmic growth of the ellipticity ratio. This mechanism allows us to close the iteration procedure leading to the local boundedness of the gradient and ultimately yields the Lipschitz estimate \eqref{estfin}.

To the best of our knowledge, this is the first result establishing local Lipschitz continuity for minimizers of variational integrals with almost linear growth under borderline Orlicz--Sobolev assumptions on the coefficients. Besides extending the recent theory developed for polynomially growing non-uniformly elliptic energies, our analysis provides further evidence that the interplay between the ellipticity ratio and the regularity of the coefficients constitutes the fundamental principle governing regularity in the non-uniformly elliptic setting.

\section{Preliminaries and notations}
\medskip
In this paper we shall denote by $C$  a
general positive constant that may vary on different occasions, even within the
same line of estimates.
Relevant dependencies  will be suitably emphasized using
parentheses or subscripts.   \\
In what follows, $B(x,r)=B_r(x)=\{y\in \R^n:\,\, |y-x|<r\}$ will denote the ball centered at $x$ of radius $r$. We shall omit the dependence on the center and on the radius when no confusion arises.
We recall the following  iteration lemma, which is a slightly different version of  the classical one (\cite[Lemma 6.1, p.191]{Giusti}) and whose
proof can be found in \cite[Lemma 2.2]{GiaPdN} 

  \begin{lemma}\label{lem:Giaq}
    For $0<R_1<R_2$, consider a bounded function
    $f:[R_1,R_2]\to[0,\infty)$ with
    \begin{equation*}
      f(r_1)\le\vartheta f(r_2)+\frac A{(r_2-r_1)^{{\sigma}}}+ \frac B{(r_2-r_1)^{{\sigma}}}\log^{\sigma}\left(1+\frac{1}{r_2-r_1}\right)+ C       
    \end{equation*}
   for all $R_1\le r_1<r_2\le R_2$, where $A,B,C$ and $\sigma$ denote nonnegative constants,
    and $\vartheta\in (0,1)$. Then there exists a constant ${c=c(\sigma, \vartheta)}>0$ such that
    \begin{equation*}
      f(R_1)\le c({\sigma, \vartheta})
     \bigg(\frac{A}{(R_2-R_1)^{{\sigma}}}+\frac{B}{(R_2-R_1)^{{\sigma}}}\log^{\sigma}\left(1+\frac{1}{R_2-R_1}\right)+C
      \bigg).
    \end{equation*}
  \end{lemma}
\subsection{Some background material on Orlicz spaces}\label{Orlicz}
Here we collect definitions and some basic properties of the {Young} functions, that will be used in what follows {({for more details} see for example \cite{KR, RR})}. 

A function {$\phi:[0,\infty)\to [0,\infty)$ is called a
Young function if and only if it is a left-continuous convex function  that vanishes at $0$ and is not constant. 
It is said to
satisfy the $\Delta_2$-condition if and only if there is a constant
$c>1$ such that $\phi(2t)\le c\,\phi(t)$ for every $t>0$. The smallest
constant $c$ with this property can be denoted by $\Delta_2(\phi)$.

The conjugate of a {Young} function~$\phi$ is defined as
\begin{align*}
  \phi^*(t) := \sup_{s \geq 0} \{ st - \phi(s) \}, \quad t \geq 0.
\end{align*}
One has that also $ \phi^*$ is  a Young function and it holds
\begin{equation}
  \label{Young-psi}
  st\le \phi(s)+ \phi^*(t)
  \qquad\mbox{for $s,t\ge0$}.
\end{equation}
The function $\phi$ is said to dominate $\psi$ near infinity  
 if there exists a constant $t_0>0$ such that $\psi (t) \leq \phi(ct)$, for $t \geq t_0$.
 We say that the functions $\phi$ and $\psi$ are equivalent globally near infinity if they dominate each other globally near infinity. The equivalence of Young functions in this sense will be denoted by $\phi \approx \psi$.
 \\
 We explicitly recall that if  $$\phi (t)\approx t \log^\beta (1+t)\quad \text{near infinity}$$ and $\beta >0$, then 
$$\phi^* (t) \approx \exp (t^{\frac{1}{\beta}})-1\quad \text{near infinity.}$$
The  Orlicz space $L^\phi(\Omega)$, generated by a Young function $\phi$, is defined as 
\begin{align*}
    \label{orlicz}
    L^\phi (\Omega)= \bigg\{u: \Omega \to \mathbb R^n: \exists \lambda >0\,\,\text{s.t.} \,\, \int_\Omega\phi\bigg(\frac{|u|}{\lambda}\bigg)\, \dx <\infty\bigg\}
\end{align*}
and it becomes a Banach space with the Luxemburg norm
\begin{equation}
\label{Luxemburgnorm}
  \|u\|_{L^\phi(\Omega)}:=\inf\bigg\{\lambda>0\colon \int_\Omega
  \phi\bigg(\frac{|u(x)|}{\lambda}\bigg)\dx\le 1\bigg\}.
\end{equation}
The Orlicz-Sobolev space $W^{1,\phi}(\Omega)$ is defined as the space of
functions $u\in L^\phi(\Omega)$ that are weakly differentiable
with $Du\in  L^\phi(\Omega,\R^n)$. It is  equipped with the
norm
$$\|u\|_{W^{1,\phi}(\Omega)}:=\|u\|_{L^{\phi}(\Omega)}+\|Du\|_{L^{\phi}(\Omega,\R^{n})}.$$ 
We shall use the common notations
$L^p\log ^\beta L (\Omega)$ and ${\rm EXP}_\alpha(\Omega)$ to denote the Orlicz spaces associated with the functions $t^p \log^\beta (1+t)$ and $\exp (t^{\alpha })-1$, respectively. We shall also use the symbol $W^{1,L\log ^\beta L}(\Omega)$ to denote the Orlicz Sobolev space associated with the functions $t \log^\beta (1+t)$.
\\
In the sequel we will use the following {equality}
\begin{equation}
\label{galquadrato}
\|\phi^2\|_{{\rm EXP}_{\alpha}(\Omega)} = \|\phi\|^2_{{\rm EXP}_{2 \alpha}(\Omega)},
\end{equation}
that easily follows by the definition of the Luxemburg norm \eqref{Luxemburgnorm}, 
For future reference,  we also observe that for every {$\beta>1$} there exists ${C(\beta)>0}$ such that 
{\begin{equation}\label{alphallogl}||f||_{L\log L(\Omega)}\le \int_{\Omega}|f|\log\left(e+|f|\right)\,\dx +C(\beta)\left(\int_{\Omega}|f|\log\left(e+|f|\right)\,\dx\right)^{\frac{1}{{\beta}}}\end{equation}}
 For an explicit proof of \eqref{galquadrato} and \eqref{alphallogl} we refer to \cite{EPdN}. 
\\
The following H\"older type inequality in Orlicz spaces will be useful in the sequel.
\begin{proposition}\label{HoldOrlicz} Let $\phi$ and $\phi^*$ be a pair of conjugate Young functions, and let $f \in L^\phi(\Omega)$ and $g \in L^{\phi^*}(\Omega)$. Then \[ \int_\Omega |f(x)g(x)|\, \dx\ \leq 2\, \|f\|_{L^\phi(\Omega)} \, \|g\|_{L^{\phi^*}(\Omega)}, \] where $\|\cdot\|_{L^\phi(\Omega)}$ and $\|\cdot\|_{L^{\phi^*}(\Omega)}$ denote the Luxemburg norms associated with $\phi$ and $\phi^*$, respectively. \end{proposition}
We shall use  the following well known  Moser Trudinger type imbedding theorem in the Orlicz Sobolev setting (\cite {cianchi_IUMJ, trudinger}).
\begin{theorem}\label{sobtru}
Let $\Omega$ be an open set in $\mathbb{R}^n$ with $|\Omega|<\infty$. Then, there exists a constant $\mathcal{C}_S=\mathcal{C}_S(n,|\Omega|)$ such that
    \begin{align}
        \label{SPnew1}
        ||w||_{{\rm EXP}_{\frac{n}{n-1}}(\Omega)} \leq \mathcal{C}_S\left( \int_\Omega |D w|^n \dx\right)^{\frac{1}{n}}
    \end{align}
    for every $w \in W^{1,n}_0(\Omega)$.
\end{theorem}
We now recall a classical embedding relating the $L^p$
 spaces to spaces of exponentially integrable functions.  
\begin{proposition}\label{normalpexp}
Let $\Omega$ be an open set in $\mathbb{R}^n$ with $|\Omega|<\infty$, let $\alpha>0$ and let $f\in \mathrm{EXP}_\alpha(\Omega)$. Then there exists a constant
$C=C(\alpha,|\Omega|)>0$ such that
\[
\|f\|_{L^p(\Omega)}
\le C\, p^{1/\alpha}\, \|f\|_{ \mathrm{EXP}_\alpha(\Omega)},
\]
for every $p\ge 1$.
\end{proposition}
\subsection{Some elementary inequalities}
In this subsection, we record some inequalities involving the logarithm, that will be useful for our aims.
We start with the following
\begin{lemma}
For every $\alpha>2$, $\mathcal{A}>0$ and $q>0$, there exists a positive constant $C=C(q,\alpha,\mathcal{A})$ such that
\begin{equation}\label{ordinelognuovo}
\Big(1+\log^2 t\Big)^q\exp\big( \mathcal{A}(1 + \log^2 t)^{\frac{1}{\alpha}}\big)
  \le \, C(q)\frac{ t}{\log t}+ {C(q,\alpha).}
\end{equation}
for every $t> 1$.
\end{lemma}
\begin{proof} We start proving that, for every $\alpha>2$, ${\mathcal{A} > 0}$ and $q>0$, there exists $\sigma=\sigma(q,\alpha,\mathcal{A})$ such that 
\begin{equation}\label{ordinedilog0}\log^{2q} t\cdot\exp(\mathcal{A}\log^{\frac{2}{\alpha}} t)\le \frac{t}{\log t}, \end{equation}
for every $t>e^\sigma$. Previous inequality  is equivalent to the following
\begin{equation}\label{ordinedilog}\log^{2q+1} t\cdot\exp(\mathcal{A}\log^{\frac{2}{\alpha}} t)\le t\qquad {\forall t > e^{{\sigma}}}\end{equation}
 that, taking the logarithm at both sides, becomes
\[
(2q+1)\log\log t+\mathcal{A}\log^{\frac{2}{\alpha}} t\le \log t.
\]
Since both the functions
$$\varphi_1(t)= \frac{(2q+1)\log\log t}{\log t}\quad\text{and}\quad \varphi_2(t)=\frac{\mathcal{A}\log^{\frac{2}{\alpha}} t}{\log t}$$
are strictly decreasing in ${(e^{\sigma},+\infty)}$, with $\sigma>e$, it holds
$$\varphi_1(t)\le \varphi_1(e^{{\sigma}})={\frac{(2q+1)\log\sigma}{\sigma}},\qquad\qquad  {\forall t > e^{{\sigma}}}$$
as well as
$$\varphi_2(t)\le\varphi_2(e^{{\sigma}})=\mathcal{A}{\sigma^{\frac{2}{\alpha}-1}}, \qquad \qquad\qquad \forall {t > {e^\sigma}.}$$
Hence
\[
\frac{(2q+1)\log\log t+\mathcal{A}\log^{\frac{2}{\alpha}} t}{\log t}\le \frac{(2q+1)\log\sigma+\mathcal{A}\sigma^{\frac{2}{\alpha}}}{\sigma }\le (c_\alpha(2q+1)+\mathcal{A})\sigma^{\frac{2}{\alpha}-1},
\]
{where we used that $\log \sigma\le c_\alpha \sigma^{\frac{2}{\alpha}}$, for $\sigma>e$.}
Hence,
 {in order to have $$(c_\alpha(2q+1)+\mathcal{A})\sigma^{\frac{2}{\alpha}-1}\le 1\quad\Longleftrightarrow\,\, \sigma^{1-\frac{2}{\alpha}}\ge (c_\alpha(2q+1)+\mathcal{A}) $$
 it suffices to choose $\sigma= e+(c_\alpha(2q+1)+\mathcal{A})^{\frac{\alpha}{\alpha-2}}$. With this choice } it follows
\[
(2q+1)\log\log t+\mathcal{A}\log^{{\frac{2}{\alpha}}} t\le\log t \qquad \forall {t > e^\sigma.}
\]
Exponentiating both sides yields \eqref{ordinedilog}, {for every $\mathcal{A} >0,$ so let write it with $\mathcal{B}$ instead of $\mathcal{A};$ from this} 
 we deduce that
for every $\sigma>e$, it holds 
\begin{equation*}
  \begin{cases}  \Big(1+\log^2t\Big)^q\exp\Big( \mathcal{B}(1 + \log^2t)^{\frac{1}{\alpha}}\Big)\le (1+\sigma^2)^q\exp\Big( \mathcal{B}(1 + \sigma^2)^{\frac{1}{\alpha}}\Big)\qquad\qquad\qquad\text{if}\,\,1<t\le e^\sigma,\\
  \\
  \Big(1+\log^2t\Big)^q\exp\Big( \mathcal{B}(1 + \log^2t)^{\frac{1}{\alpha}}\Big)\le (1 +\frac{1}{\sigma^2})^q \log^{2q} t\exp\Big( \mathcal{B}(1 +\frac{1}{\sigma^2})^{\frac{1}{\alpha}} \log^{\frac{2}{\alpha}} t\Big)\\
  \qquad\qquad\qquad\qquad\qquad\le 2^q \log^{2q} t\exp\Big( \mathcal{B}2^{\frac{1}{\alpha}} \log^{\frac{2}{\alpha}} t\Big)\qquad\qquad\qquad\qquad\qquad\qquad\text{if}\,\,t> e^\sigma.\end{cases}
\end{equation*}
Choosing $\sigma$ so that \eqref{ordinedilog0} holds with $\mathcal{A}=\mathcal{B}2^{\frac{1}{\alpha}}$, the previous inequalities yield 
\begin{equation*}
  \begin{cases}  \Big(1+\log^2t\Big)^q\exp\Big( \mathcal{A}(1 + \log^2 t)^{\frac{1}{\alpha}}\Big)\le C(q,\alpha)\qquad\qquad\qquad\qquad\text{if}\,\,t\le e^{\sigma},\\
  \\
  \Big(1+\log^2t\Big)^q\exp\Big( \mathcal{A}(1 + \log^2 t)^{\frac{1}{\alpha}}\Big)
  \le 2^q\displaystyle{\frac{ t }{\log t}}\qquad\qquad\qquad\qquad\,\,\,\text{if}\,\,t> e^\sigma.\end{cases}
\end{equation*}
{We conclude combining the two previous inequalities.}
\end{proof}
We also have the following
\begin{lemma}\label{PdN}
    For every $s, t \ge 0$ and for every $\alpha > 0$ we have that
$$st \le \sigma s \log^\alpha(1 + \beta s) + \frac{t}{\beta}\exp\left(\left(\frac{t}{\sigma}\right)^{\frac{1}{\alpha}}-1\right), $$
where $\sigma,\beta > 0$.
\end{lemma}
For the proof we refer \cite[Lemma 2.5]{PdN15}. A further observation is contained in the following inequality
\begin{equation}\label{elem}
\frac{t}{\log t}\log\left(1+\frac{t}{\log t}\right)\le \log(1+e)\cdot  t, \qquad \text{for all}\,\, t\ge e. 
\end{equation}

Remark that  { choosing $\alpha=n$ and $\beta=\frac{4(2n)^n}{n}$ and $\sigma=2^n$ in Lemma \ref{PdN} it } can also be deduced that
\begin{equation}\label{elembis}st \le  \frac{1}{4}\exp\left(t^\frac{1}{n}\right)+ C_n s \log^n(1 + C_n s) \end{equation}
for every $s,t>0$.

\section{A higher differentiability result}

Here, we recall an higher differentiability result for minimizers of \eqref{functional}, under an Orlicz-Sobolev assumption on the partial map $x\to F(x,\xi)$ that will be useful in the approximation argument.
\begin{theorem}\label{giova}
Let $u\in W^{1,{1}}_{\mathrm{loc}}(\Omega)
$
be a local minimizer of the functional at \eqref{functional}
under assumptions {{\rm (F1)--(F4)}.}
Then
\[
\frac{Du}{\sqrt{1+|Du|}}\in W^{1,2}_{\mathrm{loc}}(\Omega),
\]    
\end{theorem}
For the proof we refer to \cite[Theorem 1.1]{Giova}.

\begin{remark}\label{remark} If $u\in W^{1,\infty}_{\mathrm{loc}}(\Omega)
$
is a local minimizer of the functional at \eqref{functional}
under assumptions {{\rm (F1)--(F4)}}, then  \[u\in W^{2,2}_{\mathrm{loc}}(\Omega) . \]
Indeed, it suffices to observe that  the previous Theorem yields
\[\int_{B_r}\frac{|D^2u|^2}{(1+|Du|)}\,\dx<+\infty \qquad \text{for all}\,\, B_r\Subset\Omega\]
and so, by the assumption $u\in W^{1,\infty}_{\mathrm{loc}}(\Omega)
$
\[\int_{B_r}|D^2u|^2\,\dx=\int_{B_r}\frac{|D^2u|^2}{(1+|Du|)}(1+|Du|)\,\dx\]
\[\le ||1+|Du|||_{L^{\infty}(B_r)}\int_{B_r}\frac{|D^2u|^2}{1+|Du|}\,\dx<+\infty.\]
\end{remark}

\section{A Caccioppoli type inequality}
In this section we  prove a second order Caccioppoli type inequality under
the a priori assumption $u\in W^{1,\infty}_{\loc}(\Omega)$.
Actually, thanks to Remark \ref{remark} this assumption implies that $u\in W^{2,2}_{\loc}(\Omega)$ that allows us to  test the differentiated Euler Lagrange of the functional $\mathcal{F}(u)$, with functions  proportionals to the gradient of the minimizers and this represents the first step in a Moser type iteration argument.

\noindent To shorten the notations, from now on, we set
\begin{equation}
\label{V(Du)}
V(Du)=(1+|Du|^2)^{\frac{1}{2}}.
\end{equation}
\begin{lemma}
Let $u\in W^{1,\infty}_{\loc}(\Omega)$  be  a local minimizer of the functional $\mathcal{F}(u),$ under assumptions {{\rm (F1)--(F4).}} Then the following second order Caccioppoli type inequality
\begin{eqnarray}\label{Caccioppoli}
 && \int_{\Omega} \eta^2 V(Du)^{\gamma-1}  |D^2 u|^2 \, \dx \cr\cr
&\le & C (1+\gamma^2)  \, \int_{\Omega} \eta^2  h^2(x) V(Du)^{\gamma + 1} (1 + \log^2 V(Du) ) \, \dx\cr\cr
&&+ C \int_{\Omega} |D \eta|^2 \,  V(Du)^{\gamma + 1} \,(1+\log V(Du))   \, \dx
\end{eqnarray}
holds true with a constant $C= C(\lambda,\Lambda, n)>0$, for every {$\gamma \ge \, 0$} and for every $\eta\in C^1_0(\Omega).$
\end{lemma}
\begin{proof}
{By the assumption 
$
u \in  W^{1,\infty}_{\loc}(\Omega)$ and by Remark \ref{remark}, we have that $u\in  W^{2,2}_{\loc}(\Omega)$ and so the following system
\begin{equation}
\label{second_variation}
\int_{\Omega} \left (\sum_{i,j=1}^n  F_{\xi_i \xi_j}(x, Du)u_{x_j x_s}D_{x_i} \varphi   + \sum_{i=1}^n F_{\xi_i x_s}(x, Du) D_{x_i} \varphi \right )\, \dx = 0,
\end{equation}
holds for all $s = 1, \dots, n$ and for all $\varphi \in W^{1,2}_0(\Omega)$. For $\eta \in \mathcal{C}^1_0(\Omega)$ and $\gamma \ge 0$, let us consider the function $\varphi=\eta^2
(1+|Du|^2)^{\frac{\gamma}{2}} u_{x_s}={\eta^2
V(Du)^{\gamma} u_{x_s}}$,   so that
\begin{eqnarray*}
 D_{x_i} \varphi &=&2 \eta \eta_{x_i} V(Du)^{\gamma} u_{x_s} + \gamma \eta^2  V(Du)^{\gamma - 2}|Du| D_{x_i}(|Du|) u_{x_s} \\
 &&+ \eta^2 V(Du)^{\gamma} u_{x_s x_i}. 
\end{eqnarray*} 
Using $\varphi$ as test function in \eqref{second_variation}, we get
\begin{eqnarray*}
0 &=& 2\int_{\Omega} \eta\sum_{i,j=1}^n F_{\xi_i \xi_j}(x, Du) u_{x_j x_s}  \eta_{x_i} V(Du)^{\gamma} u_{x_s} \, \dx\\
&& +  \int_{\Omega} \eta^2\sum_{i,j=1}^n F_{\xi_i \xi_j}(x, Du) u_{x_j x_s} V(Du)^{\gamma} u_{x_s x_i} \, \dx\\
&& + \gamma \int_{\Omega}\eta^2 \sum_{i,j=1}^n F_{\xi_i \xi_j}(x, Du) u_{x_j x_s}   V(Du)^{\gamma - 2}|Du| D_{x_i}(|Du|) u_{x_s} \, \dx\\
&& +2 \int_{\Omega} \eta\sum_{i=1}^n F_{\xi_i x_s}(x, Du)  \eta_{x_i}  V(Du)^{\gamma} u_{x_s} \, \dx\\
&& +  \int_{\Omega} \eta^2\sum_{i=1}^n F_{\xi_i x_s }(x, Du) V(Du)^{\gamma} u_{x_s x_i} \, \dx\\
&& + \gamma \int_{\Omega} \eta^2\sum_{i=1}^n F_{\xi_i x_s}(x, Du)  V(Du)^{\gamma - 2} |Du| D_{x_i}(|Du|) u_{x_s} \, \dx\\
&=:& I_1+I_2+I_3+I_4+I_5+I_6.
\end{eqnarray*}
Summing all terms in the previous  equation   for $s=1,\dots,n$ , and  denoting once more with ${I_1-I_6}$ the corresponding integrals, we obtain
\begin{equation}\label{start}
	I_2+I_3\le |{I}_1|+|{I}_4|+|{I}_5|+|{I}_6|.
\end{equation}
Recalling that $D_{x_j}(|Du|)|Du|={\sum_{s=1}^n u_{x_j x_s} u_{x_s}}$, for the term $I_3$  we have
\begin{eqnarray*}\label{I3}
{I}_3 &=& \gamma\int_{\Omega} \eta^2 \sum_{i,j,s =1}^n F_{\xi_i \xi_j}(x, Du) u_{x_j x_s} \left [   V(Du)^{\gamma - 2}  D_{x_i}(|Du|)|Du| \right ] u_{x_s} \, \dx\cr\cr
&\ge&  \gamma \, \int_{\Omega} \eta^2 V(Du)^{\gamma - 2}|Du|^2 \sum_{i,j =1}^n F_{\xi_i \xi_j}(x, Du) D_{x_i}(|Du|) D_{x_j}(|Du|) \, \dx\cr\cr
&\ge&  \gamma \,  {\lambda} \int_{\Omega} \eta^2 V(Du)^{\gamma - 2} |Du|^2  \frac{|D(|Du|)|^2}{  (1 + |Du|) } \, \dx \ge 0,
\end{eqnarray*}
where, in the last line we used the left inequality in {(F2).}
Hence, estimate \eqref{start} implies
\begin{eqnarray}\label{ristart}
	I_2\le |{I}_1|+|{I}_4|+|{I}_5|+|{I}_6|.
\end{eqnarray}
Cauchy-Schwarz and Young's inequalities yield
\begin{eqnarray}\label{I1}
|{I}_1| &=& 2\left |\int_{\Omega}  \eta  V(Du)^{\gamma}  \sum_{i,j,s =1}^n F_{\xi_i \xi_j}(x, Du) u_{x_j x_s} \eta_{x_i} u_{x_s} \, \dx\right |\cr\cr
&\le& 2\int_{\Omega}  \eta V(Du)^{\gamma}  \cr\cr
&& \times \left \{ \sum_{i,j,s =1}^n F_{\xi_i \xi_j}(x, Du) \eta_{x_i} \eta_{x_j} u_{x_s}^2\right \}^{1/2} \, \left \{ \sum_{i,j,s =1}^n F_{\xi_i \xi_j}(x, Du) u_{x_s x_i} \, u_{x_s x_j}\right \}^{1/2} \, \dx \cr\cr
&\le &  C   \int_{\Omega}  V(Du)^{\gamma} \sum_{i,j,s =1}^n F_{\xi_i \xi_j}(x, Du) \eta_{x_i} \eta_{x_j} u_{x_s}^2 \, \dx\cr\cr
&&  + \frac{1}{2} \int_{\Omega} \eta^2  V(Du)^{\gamma}
\sum_{i,j,s =1}^n F_{\xi_i \xi_j}(x, Du) u_{x_s x_i} \, u_{x_s x_j}  \, \dx \cr\cr
&\le& {C(\Lambda)} \, \int_{\Omega} |D \eta|^2{\frac{\log(1+|Du|)}{|Du|}}|Du|^2 \,V(Du)^{\gamma} \, \dx \cr\cr
&& + \frac{1}{2} \int_{\Omega} \eta^2  V(Du)^{\gamma} \, \sum_{i,j,s =1}^n F_{\xi_i \xi_j}(x, Du) u_{x_j x_s} u_{x_i x_s} \, \dx\cr\cr 
&\le& {C(\Lambda)} \, \int_{\Omega} |D \eta|^2 \, V(Du)^{\gamma + 1} (1 + \log(V(Du)))\, \dx \cr\cr
&& + \frac{1}{2} \int_{\Omega} \eta^2  V(Du)^{\gamma}  \, \sum_{i,j,s =1}^n F_{\xi_i \xi_j}(x, Du) u_{x_j x_s} u_{x_i x_s} \, \dx,
\end{eqnarray}
where, in the second last inequality, we used the right inequality in assumption  {(F2)} and where, from now on, we will use sistematically the elementary inequality
\[
\log(1 + |Du|) \le \, \log \Big((2 (1 + |Du|^2))^{1/2}\Big) = \log (\sqrt{2} V(Du)) \le \, C (1 + \log V(Du))
\]
and therefore also
\[
\log^2(1 + |Du|) \le \, C (1 + \log^2 V(Du)).
\]
Using again Young's inequality together with {(F3)} we derive that
\begin{eqnarray}\label{I4}
|{I}_4| &=& 2 \left | \int_{\Omega} \eta  \, V(Du)^{\gamma} \sum_{i,s =1}^n F_{\xi_i x_s}(x, Du) \eta_{x_i} u_{x_s} \, \dx \right |\cr\cr
&\le&   C \,  \int_{\Omega} \eta |D\eta| |Du| \, h(x) \, V(Du)^{\gamma}{\log(1+|Du|) }  \, \dx \cr\cr
&\le& C \,  \int_{\Omega} \eta |D\eta|  \, h(x) \, V(Du)^{\gamma + 1}{\log(1+|Du|) }  \, \dx \cr\cr
&\le& C\int_{\Omega} |D\eta|^2  V(Du)^{\gamma + 1}\,\dx\\\nonumber
&&+C\int_{\Omega} \eta^2 h^2(x) V(Du)^{\gamma + 1} \log^2(1+|Du|) \,\dx \cr\cr
&\le& C\int_{\Omega} |D\eta|^2  V(Du)^{\gamma + 1}\,\dx\\\nonumber
&&+C\int_{\Omega} \eta^2 h^2(x) V(Du)^{\gamma + 1} (1 + \log^2V(Du)) \,\dx,
\end{eqnarray}
and also, for a parameter $\varepsilon>0$ that will be chosen later,
\begin{eqnarray}\label{I5}
|{I}_5| &=& \left| \int_{\Omega} \eta^2 V(Du)^{\gamma} \sum_{i,s=1}^n F_{\xi_i x_s }(x, Du)  u_{x_s x_i} \, \dx \right| \cr\cr
&\le&\,  \, \int_{\Omega} \eta^2 h(x) V(Du)^{\gamma} {\log(1+|Du|)}  |D^2 u| \, \dx \cr\cr
&=& \int_{\Omega} \left[\eta^2 V(Du)^{\gamma}  {\frac{|D^2 u|^2}{(1+|Du|)} }\right]^{\frac{1}{2}} \left[\eta^2 h^2(x) V(Du)^{\gamma}{(1+|Du|)\log^2(1+|Du|)}\right]^{\frac{1}{2}} \, \dx \cr\cr
&\le& \varepsilon \int_{\Omega}  \eta^2 V(Du)^{\gamma}  {\frac{|D^2 u|^2}{(1+|Du|)} }\, \dx\\\nonumber&& + C_\varepsilon  \int_{\Omega} \eta^2 h^2(x) V(Du)^{\gamma + 1}{\log^2(1+|Du|)}  \, \dx \cr\cr
&\le& \varepsilon \int_{\Omega}  \eta^2 V(Du)^{\gamma} {\frac{|D^2 u|^2}{(1+|Du|)} \, \dx}\\\nonumber&& + C_\varepsilon  \int_{\Omega} \eta^2 h^2(x) V(Du)^{\gamma + 1} (1 + \log^2 V(Du))  \, \dx.
\end{eqnarray}
Similar arguments give
\begin{eqnarray}\label{I6}
|{I}_6| &=& \gamma \,  \left | \int_{\Omega} \sum_{i,s=1}^n F_{\xi_i x_s}(x, Du) \eta^2  V(Du)^{\gamma - 2}|Du| D_{x_i}(|Du|) u_{x_s} \, \dx \right |\cr\cr
&\le&  \, \gamma \, \int_{\Omega} \eta^2 V(Du)^{\gamma - 2} {|Du|}\left |\sum_{i,s=1}^n F_{\xi_i x_s}(x, Du)  D_{x_i}(|Du|) u_{x_s} \, \dx \right |\cr\cr
&\le&  \, \gamma \, \int_{\Omega} \eta^2  h(x) V(Du)^{\gamma - 2}{\log(1+|Du|)} {|Du|}\left | \sum_{i,s=1}^n D_{x_i}(|Du|) u_{x_s} \, \dx \right |\cr\cr
&\le& C \, \gamma \, \int_{\Omega} \, \eta^2 V(Du)^{\gamma}{\log(1+|Du|)}|D^2 u| \, h(x)\, \dx\cr\cr
&\le& \varepsilon \int_{\Omega} \eta^2 \frac{|D^2u|^2 }{{(1+|Du|)}} V(Du)^{\gamma}\, \dx\\\nonumber
&&+ C_\varepsilon \gamma^2 \,  \, \int_{\Omega} \eta^2  h^2(x) V(Du)^{\gamma + 1} (1 + \log^2 V(Du)) \, \dx,
\end{eqnarray}
where  the constants $C$ and $C_{\varepsilon}$ are independent of $\gamma$.}
Inserting  \eqref{I1}, \eqref{I4}, \eqref{I5} and \eqref{I6} in \eqref{ristart} we obtain
\begin{eqnarray*}
&& \int_{\Omega} \eta^2 V(Du)^{\gamma}{\sum_{i,j,s=1}^n} F_{\xi_i \xi_j}(x, Du) u_{x_j x_s}  u_{x_s x_i} \, \dx\cr\cr
&\le&  \frac{1}{2} \int_{\Omega} \eta^2  V(Du)^{\gamma} \, \sum_{i,j,s =1}^n F_{\xi_i \xi_j}(x, Du) u_{x_j x_s} u_{x_i x_s} \, \dx\cr\cr
&&+2\varepsilon \int_{\Omega} \eta^2 V(Du)^{\gamma}  \frac{|D^2u|^2 }{{(1+|Du|)}} \, \dx \cr\cr
&&+ C_\varepsilon (1+\gamma^2)  \, \int_{\Omega}  \eta^2  h^2(x) V(Du)^{\gamma + 1} (1 + \log^2 V(Du))\, \dx\cr\cr
&&+C_\varepsilon \int_{\Omega} |D \eta|^2 \, V(Du)^{\gamma + 1}(1+ \log V(Du))  \, \dx.
\end{eqnarray*}
Reabsorbing the first integral in the right hand side by the left hand side, we obtain
\begin{eqnarray*}
&& \frac{1}{2}\int_{\Omega}\eta^2 V(Du)^{\gamma} {\sum_{i,j,s=1}^n} F_{\xi_i \xi_j}(x, Du) u_{x_j x_s}  u_{x_s x_i} \, \dx\cr\cr
 &\le &2\varepsilon \int_{\Omega} \eta^2 V(Du)^{\gamma} \frac{|D^2u|^2 }{{(1+|Du|)}} \, \dx \cr\cr
&&+ C_\varepsilon (1+\gamma^2)  \, \int_{\Omega}  \eta^2  h^2(x) V(Du)^{\gamma + 1} (1 + \log^2 V(Du)) \, dx\cr\cr
&&+C_\varepsilon \int_{\Omega} |D \eta|^2 \, V(Du)^{\gamma + 1} (1+\log V(Du))   \, \dx.
\end{eqnarray*}
Using the left inequality in assumption {(F2)} in the left hand side of previous estimate, we deduce
\begin{eqnarray*}
&& \frac{{\lambda}}{2}\int_{\Omega} \eta^2 V(Du)^{\gamma}{\frac{|D^2 u|^2}{(1+|Du|)}} \, \dx \cr\cr
 &\le &2\varepsilon \int_{\Omega} \eta^2 V(Du)^{\gamma} \frac{|D^2u|^2 }{{(1+|Du|)}} \, \dx \cr\cr
&&+ C_\varepsilon (1+\gamma^2)  \, \int_{\Omega}  \eta^2  h^2(x) V(Du)^{\gamma + 1} (1 + \log^2 V(Du)) \, \dx\cr\cr
&&+C_\varepsilon \int_{\Omega} |D \eta|^2 \, V(Du)^{\gamma + 1} (1+ \log V(Du))  \, \dx.
\end{eqnarray*}
Choosing, for instance, $\varepsilon=\frac{{\lambda}}{8}$, we can reabsorb the first integral in the right hand side by the left hand side thus getting
\begin{eqnarray*}
&& \int_{\Omega} \eta^2 V(Du)^{\gamma - 1} |D^2 u|^2 \, \dx \cr\cr
 &\le & C (1+\gamma^2)  \, \int_{\Omega}  \eta^2  h^2(x) V(Du)^{\gamma + 1} (1 + \log^2 V(Du)) \, \dx\cr\cr
&&+C \int_{\Omega} |D \eta|^2 \,  V(Du)^{\gamma + 1}(1+\log V(Du))  \, \dx,
\end{eqnarray*}
with a constant $C=C(\lambda, \Lambda, n)>0$, and this concludes the proof.
\end{proof}

\medskip

\section{The a priori estimate}

The aim of this section is to prove a uniform a priori estimate for the Lipschitz norm of the minimizers of the functional at \eqref{functional}, by means of a Moser type iteration argument. The proof will be divided in two parts, distiguishing the cases $n\ge 3$ and $n=2$, that reflect the use of two different Sobolev embedding theorem. More precisely, we are going to prove the following

\begin{theorem}\label{thmapriori}
Let $u\in W^{1,\infty}_{\loc}(\Omega)$,  be  a local minimizer of the functional {at \eqref{functional}}, under the assumptions {{\rm (F1)--(F4)}}. Then the following estimate \begin{eqnarray}\label{apriori}
	||V(Du)||_{L^{\infty}(B_{\rho_0})} &\le&   C\int_{B_{  R_0}}V(Du)
\, \dx\cdot \exp\Bigg({C}\int_{B_{ R_0}}V(Du)\, \dx\Bigg)\cr\cr
    &&+\left(\frac{C}{(R_0- \rho_0 )^{{\theta}}}\int_{B_{ R_0}}V(Du)\, \dx\right)\log^{{\theta}}\left(1+\frac{{C}}{(  R_0- \rho_0 )^{{\theta}}}\int_{B_{  R_0}}V(Du)\, \dx\right) \cr\cr
 &&+\frac{C}{(  R_0- \rho_0 )^{{\theta}}}\int_{B_{ R_0}}V(Du)\, \dx+C\int_{B_{ R_0}}V(Du)\, \dx
\end{eqnarray}
holds true for every pair of concentric balls $B_{\rho_0}\subset B_{R_0}\Subset\Omega$ with a constant $C=C(n,\alpha,\lambda,\Lambda, h)$, and where  $V(Du)$ has been defined at \eqref{V(Du)} {and where \[
\theta =
\begin{cases}
n, & n\ge3,\\[2mm]
2q, & n=2,
\end{cases}
\]
for any $q < 1.$}
\end{theorem}

{\subsection{Proof of Theorem \ref{thmapriori} in the case $n \ge 3$}}

\begin{proof}
Let us fix a ball $B_{R_0}\Subset \Omega$ and radii $\rho_{0}<\bar\rho<\rho<R<\bar R<R_{0}.$
 We choose in estimate \eqref{Caccioppoli} a cut-off function $\eta\in C^1_0(B_{R})$ such that $\eta=1$ on $B_{\rho}$ and $|D\eta|\le \frac{C}{R-\rho}$, {thus getting}
\begin{eqnarray}\label{estprimadellasob}
&& \int_{B_{R}} \eta^2 V(Du)^{\gamma - 1} |D^2 u|^2 \, \dx \cr\cr
&\le & C (1+\gamma^2)  \, \int_{B_{R}}\eta^2  h^2(x) V(Du)^{\gamma + 1} (1 + \log^2 V(Du)) \, \dx\cr\cr
&&\quad+  \frac{C}{(R-\rho)^2}  \, \int_{B_{R}} V(Du)^{\gamma + 1}(1+\log V(Du))   \, \dx.\label{primastima}
\end{eqnarray}
The Sobolev  inequality yields
\begin{eqnarray*}
&& \left ( \int_{B_{R}}  \eta^{2^*} V(Du)^{(\gamma + 1) \frac{2^*}{2}}\, \dx\right )^{\frac{2}{2^*}}\\
& \le& \, C \, \int_{B_{R}} \Big|D (\eta V(Du)^{\frac{\gamma + 1}{2}})\Big|^2 \, \dx 	\\
&\le&C(1+\gamma^2)\int_{B_{R}} \eta^2 V(Du)^{\gamma - 1} |D^2 u|^2\,\dx+C \int_{B_{R}} |D \eta|^2 \,  V(Du)^{\gamma + 1}  \, \dx,
\end{eqnarray*}
where $2^*$  denotes as usual the Sobolev conjugate exponent of $2$, i.e., since $n\ge 3$, $2^*=\frac{2n}{n-2}$.
	Using estimate \eqref{primastima} to control the first integral in the right hand side of previous inequality, we obtain
\begin{eqnarray}\label{stimapreite}
&&\left ( \int_{B_{R}}  \eta^{2^*} V(Du)^{(\gamma + 1)\frac{2^*}{2}}\, \dx\right )^{\frac{2}{2^*}} 	\cr\cr
&\le& C(1+\gamma^2)^2 \int_{B_{R}}\eta^2  h^2(x)V(Du)^{\gamma + 1} (1 + \log^2 V(Du)) \, \dx \cr\cr
&&+\frac{C(1+\gamma^2)}{(R - \rho)^{2}} \int_{B_{R}}  \, V(Du)^{\gamma + 1}(1+\log V(Du)) \, \dx \cr\cr
&=&  C(1+\gamma^2)^2 \int_{h >\varepsilon}\eta^2
\frac{h^2}{{\log^{2\alpha }(1+h)}} {\log^{2\alpha }(1+h)} V(Du)^{ \gamma+1} (1 + \log^2V(Du))  \, \dx \cr \cr 
&& + C \, \varepsilon^2(1+\gamma^2)^2 \int_{h < \varepsilon}\eta^2 V(Du)^{ \gamma+1} (1 + \log^2 V(Du)) \, \dx\cr\cr
&&+\frac{C (1+\gamma^2)}{(R - \rho)^{2}} \int_{B_{R}}  \,  V(Du)^{ \gamma+1} (1+\log V(Du)) \, \dx \cr \cr 
&\le& C\varepsilon^2(1+\gamma^2)^2 (1 + \log^2\|{V(Du)}\|_{{L^\infty(B_{\bar R})}} )\int_{B_R} V(Du)^{\gamma+1} \, \dx\cr\cr && +\frac{C  (1+\gamma^2)(1 + \log\|{V(Du)}\|_{{L^\infty(B_{\bar R})}} )}{(R - \rho)^{2}} \int_{B_{R}}  \,  V(Du)^{ \gamma+1}  \, \dx\cr \cr 
&& + C\frac{(1+\gamma^2)^2}{{\log^{2\alpha}(1+\varepsilon)}} (1 + \log^2\|{V(Du)}\|_{{L^\infty(B_{\bar R})}} ) \int_{B_R}\eta^2 h^2 \log^{2\alpha }(1+h) V(Du)^{  \gamma+1} \,\dx\cr\cr
&=:& I + II + III,
\end{eqnarray}
where we used {the a priori assumption $u\in W^{1,\infty}_{\loc}(\Omega)$ and }the monotonicity of the function $\log^\alpha(1+t),$  where $\alpha$ is the exponent appearing in assumption {(F4)} and $\varepsilon>0$ will be chosen later.
We estimate the integral $III$ by H\"older's inequality as follows
\begin{eqnarray*}
III &\le& C\frac{(1+\gamma^2)^2}{{\log^{2\alpha}(1+\varepsilon)}} (1 + \log^2\|{V(Du)}\|_{{L^\infty(B_{\bar R})}} )\cr\cr&&\qquad\cdot \left (\int_{B_R}h^n {\log^{\alpha n}(1+h)} \, \dx \right )^{2/n} \left (\int_{B_R}\eta^{2^*}V(Du)^{\frac{( \gamma+1) n}{n-2}} \, \dx \right )^{\frac{n-2}{n}},
\end{eqnarray*}
{which is legitimate by virtue of assumption {(F4).}}
Setting
\[
H_n := \left (\int_{B_{R_0}} h^n {\log^{\alpha n}(1+h)} \, \dx \right )^{1/n},
\]
 we have that
\begin{eqnarray*}
III &\le& C\frac{H_n^2 (1+\gamma^2)^2}{{\log^{2\alpha}(1+\varepsilon)}}  \Big(1 + \log^2\|{V(Du)}\|_{{L^\infty(B_{\bar R})}} \Big)\left (\int_{B_R}\eta^{2^*} V(Du)^{\frac{( \gamma+1) n}{n-2}} \, \dx \right )^{\frac{n-2}{n}}.
\end{eqnarray*}
At this point we choose $\varepsilon$ so that
\[
\frac{C H_n^2(1+\gamma^2)^2  \Big(1 + \log^2\|{V(Du)}\|_{{L^\infty(B_{\bar R})}} \Big)}{\log^{2\alpha}(1+\varepsilon)}=\frac{1}{4}
\]
i.e.
\[
\log^{2\alpha}(1+\varepsilon)=4C H_n^2(1+\gamma^2)^2 \Big(1 + \log^2\|{V(Du)}\|_{{L^\infty(B_{\bar R})}} \Big).
\]
This yields
\[
\log(1+\varepsilon)=2^{\frac{1}{\alpha}} C^{\frac{1}{2\alpha}} H_n^{\frac{1}{\alpha}}(1+\gamma^2)^{\frac{1}{\alpha}}  \Big(1 + \log^2 \|{V(Du)}\|_{{L^\infty(B_{\bar R})}}\Big)^{\frac{1}{2\alpha}}
\]
$$\quad \Longrightarrow \quad \varepsilon=\exp\Big( 2^{\frac{1}{\alpha}} C^{\frac{1}{2\alpha}} H_n^{\frac{1}{\alpha}}(1+\gamma^2)^{\frac{1}{\alpha}}  \Big(1 + \log^2 \|{V(Du)}\|_{{L^\infty(B_{\bar R})}}\Big)^{\frac{1}{2\alpha}}\Big)-1$$
and, with this choice, {shortening $\|V(Du)\|_{{L^\infty(B_{\bar R})}}$ with $\|V(Du)\|_{\infty}$}, estimate \eqref{stimapreite} becomes
\begin{eqnarray*}
&&\left (\int_{B_{{R}}} \eta^{2^*}V(Du)^{\frac{( \gamma+1) n}{n-2}} \,\dx \right )^{\frac{n-2}{n}}\cr\cr
&\le& \frac{1}{4} \left (\int_{B_R}\eta^{2^*} V(Du)^{\frac{( \gamma+1) n}{n-2}} \, \dx \right )^{\frac{n-2}{n}}+\frac{C (1+\gamma^2)(1 + \log\|{V(Du)}\|_{{\infty}} }{(R - \rho)^{2}} \int_{B_{R}}  \,  V(Du)^{ \gamma+1}  \, \dx\cr\cr
&&\quad +C\bigg[(1+\gamma^2)^2 (1 + \log^2\|{V(Du)}\|_{\infty}) \\
&& \quad \cdot \exp\Big( 2^{\frac{1}{\alpha}+1} C^{\frac{1}{2\alpha}} H_n^{\frac{1}{\alpha}}(1+\gamma^2)^{\frac{1}{\alpha}}  (1 + \log^2 \|{V(Du)}\|_{\infty})^{\frac{1}{2\alpha}}\Big)\bigg] \int_{B_R} V(Du)^{ \gamma+1} \, \dx.
\end{eqnarray*}
Reabsorbing the first integral in the right hand side by the left hand side, and recalling that $\eta=1$ on $B_\rho$, we finally have
\begin{eqnarray}\label{baseiter}
&&\left (\int_{B_{\rho}} V(Du)^{\frac{( \gamma+1) n}{n-2}} \,\dx \right )^{\frac{n-2}{n}} \cr\cr
&\le& \frac{C (1+\gamma^2)(1 + \log\|{V(Du)}\|_{{\infty}} )}{(R - \rho)^{2}} \int_{B_{R}}  \,  V(Du)^{ \gamma+1}  \, \dx\cr\cr
&&\quad +C\bigg[(1+\gamma^2)^2 (1 + \log^2\|{V(Du)}\|_{\infty}) \cr\cr
&& \quad \cdot \exp\Big( 2^{\frac{1}{\alpha}+1} C^{\frac{1}{2\alpha}} H_n^{\frac{1}{\alpha}}(1+\gamma^2)^{\frac{1}{\alpha}}  (1 + \log^2 \|{V(Du)}\|_{\infty})^{\frac{1}{2\alpha}}\Big)\bigg] \int_{B_R} V(Du)^{ \gamma+1} \, \dx\cr\cr
&\le& \frac{C (1+\gamma^2)(1+\log  \|{V(Du)}\|_{\infty}) }{(R - \rho)^{2}} \int_{B_{R}}  \,  V(Du)^{ \gamma+1} \, \dx\cr\cr
&& +C\bigg[(1+\gamma^2)^2 (1 + \log^2\|{V(Du)}\|_{\infty})\exp\Big( \tilde C(1+\gamma^2)^{\frac{1}{\alpha}}  (1 + \log^2 \|{V(Du)}\|_{\infty})^{\frac{1}{2\alpha}}\Big)\bigg] \cr\cr
&& \quad \cdot  \int_{B_R} V(Du)^{ \gamma+1} \, \dx,
\end{eqnarray}
where we  set $\tilde C=2^{\frac{1}{\alpha}+1} C^{\frac{1}{2\alpha}} H_n^{\frac{1}{\alpha}}$. 
Using Young's inequality to deduce that $\exp(ab)\le \exp(a^2)\exp(b^2)$, from previous estimate we obtain
\begin{eqnarray}\label{baseiter}
&&\left (\int_{B_{\rho}} V(Du)^{\frac{( \gamma+1) n}{n-2}} \, \dx \right )^{\frac{n-2}{n}} \cr\cr
&\le& \frac{C (1+\gamma^2)(1+\log  \|{V(Du)}\|_{\infty}) }{(R - \rho)^{2}} \int_{B_{R}}  \,  V(Du)^{ \gamma+1} \, \dx\cr\cr
&& +C\bigg[(1+\gamma^2)^2\exp\Big( {\tilde C^2}(1+\gamma^2)^{\frac{2}{\alpha}}\Big) (1 + \log^2\|{V(Du)}\|_{\infty})\exp\Big(   (1 + \log^2 \|{V(Du)}\|_{\infty})^{\frac{1}{\alpha}}\Big)\bigg] \cr\cr
&& \qquad\qquad \cdot  \int_{B_R} V(Du)^{ \gamma+1} \, \dx.
\end{eqnarray}
Note that previous estimate holds true for every $\rho<R$ such that $\bar\rho\le  \rho< R\le \bar R$ and,  by virtue of the a priori assumption $u\in W^{1,\infty}_{\mathrm{loc}}(\Omega)$, for every $\gamma\ge 0$. Hence, defining the decreasing sequence of radii by setting
\begin{equation}\label{raggi}
  \rho_i=\bar \rho+\frac{ \bar R- \bar\rho}{2^i}  
\end{equation}
and the following increasing sequence of exponents
\begin{equation}\label{esponenti} p_0=1\qquad{p_{i}}={p_{i-1}}\frac{2^*}{2}=\left(\frac{2^*}{2}\right)^{i},\end{equation}
 we may iterate estimate \eqref{baseiter}  on the concentric balls $B_{\rho_i}$  with exponents $p_i$,
thus obtaining
\begin{eqnarray}\label{stimapreitepassoin}
&&\left ( \int_{B_{\rho_{i+1}}}  V(Du)^{p_{i+1}}\, \dx\right )^{\frac{1}{p_{i+1}}} 	\cr\cr
&\le& \displaystyle{\prod_{j=0}^{i}}\Bigg[ Cp_j^{4}\exp\Big({\tilde C^2}p_j^{\frac{4}{\alpha}}\Big)\Big(1+\log^2\|V(Du)\|_{\infty}\Big)\exp\Big(   (1 + \log^2 \|{V(Du)}\|_{\infty})^{\frac{1}{\alpha}}\Big)\cr\cr
&&\qquad +\frac{Cp_j^2}{(\rho_i-\rho_{i+1})^2}\Big(1+\log^2\|V(Du)\|_{\infty}\Big)\Bigg]^{\frac{1}{p_j}}
\cdot\int_{B_{\bar R}} V(Du) \, \dx\cr\cr
&=& \displaystyle{\prod_{j=0}^{i}}\Bigg[ Cp_j^{4}\exp\Big({\tilde C^2}p_j^{\frac{4}{\alpha}}\Big)\Big(1+\log^2\|V(Du)\|_{\infty}\Big)\exp\Big((1 + \log^2 \|{V(Du)}\|_{\infty})^{\frac{1}{\alpha}}\Big)\cr\cr
&&\qquad +\frac{C4^{j+1}p_j^2}{(\bar R-\bar\rho)^2}\Big(1+\log^2\|V(Du)\|_{\infty}\Big)\Bigg]^{\frac{1}{p_j}}
\cdot\int_{B_{\bar R}} V(Du) \, \dx\cr\cr
&\le& \displaystyle{\prod_{j=0}^{i}}\Bigg[ Cp_j^{4}\exp\Big({\tilde C^2}p_j^{\frac{4}{\alpha}}\Big)+C4^{j+1}p_j^2\Bigg]^{\frac{1}{p_j}}\Bigg[ \Big(1+\log^2\|V(Du)\|_{\infty}\Big)\exp\Big(  (1 + \log^2 \|{V(Du)}\|_{\infty})^{\frac{1}{\alpha}}\Big)\cr\cr
&&\qquad +\frac{\Big(1+\log^2\|V(Du)\|_{\infty}\Big)}{(\bar R-\bar\rho)^2}\Bigg]^{\frac{1}{p_j}}
\cdot\int_{B_{\bar R}} V(Du) \, \dx=: \Theta_i \cdot\int_{B_{\bar R}} V(Du) \, \dx. 
\end{eqnarray}
Since
$$\Theta_i=\exp\left(\sum_{j=0}^i \frac{1}{p_j}\log\left(Cp_j^{4}\exp\Big( {\tilde C^2}p_j^{\frac{4}{\alpha}}\Big)+C4^{j+1}p_j^2\right)\right)$$
$$\cdot\exp\Bigg[\sum_{j=0}^{i} \frac{1}{p_j}\log\left(\Big(1+\log^2\|V(Du)\|_{\infty}\Big)\exp\Big(   (1 + \log^2 \|{V(Du)}\|_{\infty})^{\frac{1}{\alpha}}\Big)+\frac{\Big(1+\log^2\|V(Du)\|_{\infty}\Big)}{(\bar R-\bar\rho)^2}\right)\Bigg]$$
$$\le \exp\left(\sum_{j=0}^\infty \frac{1}{p_j}\log\left(2Cp_j^{4}\right)\right)\cdot \exp\left(\sum_{j=0}^\infty \frac{1}{p_j}\log\left(2C4^{j+1}p_j^2\right)\right)\cdot \exp\left(\sum_{j=0}^\infty \frac{1}{p_j}\log\exp\Big({\tilde C^2}p_j^{\frac{4}{\alpha}}\Big)\right)$$
$$\cdot\exp\Bigg[\sum_{j=0}^{\infty} \frac{1}{p_j}\log\left(\Big(1+\log^2\|V(Du)\|_{\infty}\Big)\exp\Big(   (1 + \log^2 \|{V(Du)}\|_{\infty})^{\frac{1}{\alpha}}\Big)+\frac{\Big(1+\log^2\|V(Du)\|_{\infty}\Big)}{(\bar R-\bar\rho)^2}\right)\Bigg]$$
$$=:\Theta_1\cdot \Theta_2\cdot \Theta_3\cdot \Theta_4,$$
where we used the elementary inequality
$\log(a+b)\le \log (2a)+\log(2 b)$, for $a,b>1$.
Note that, by the definition of $p_j$ and since $\sum_{j=0}^{{\infty}}\frac{1}{p_j}= \frac{n}{ 2} $ 
we have that
$$\Theta_1\cdot\Theta_2\le C(n,\lambda,\Lambda, H_n).$$
For the estimate of $\Theta_3$, we use the assumption $\alpha>4$ to deduce that
$$\Theta_3=\exp\left(\sum_{j=0}^\infty  \tilde Cp_j^{\frac{4}{\alpha}-1}\right)\le C(n,\alpha,\lambda,\Lambda, H_n).$$
Moreover
$$ \Theta_4\le \exp\Bigg[\frac{n}{2}\log\left(\Big(1+\log^2\|V(Du)\|_{\infty}\Big)\exp\Big(   (1 + \log^2 \|{V(Du)}\|_{\infty})^{\frac{1}{\alpha}}\Big)+\frac{\Big(1+\log^2\|V(Du)\|_{\infty}\Big)}{(\bar R-\bar\rho)^2}\right)\Bigg]$$
$$\le c(n) \Big(1+\log^2\|V(Du)\|_{\infty}\Big)^\frac{n}{2}\exp\Big( \frac{n}{2}  (1 + \log^2 \|{V(Du)}\|_{\infty})^{\frac{1}{\alpha}}\Big)+c(n)\frac{\Big(1+\log^2\|V(Du)\|_{\infty}\Big)^\frac{n}{2}}{(\bar R-\bar\rho)^n}.$$
Combining all the estimates above, we conclude that
\begin{eqnarray*}
 \Theta_i&\le& \mathfrak{C}\left(\Big(1+\log^2\|V(Du)\|_{\infty}\Big)^\frac{n}{2}\exp\Big( \frac{n}{2}  (1 + \log^2 \|{V(Du)}\|_{\infty})^{\frac{1}{\alpha}}\Big)+\frac{\Big(1+\log^2\|V(Du)\|_{\infty}\Big)^\frac{n}{2}}{(\bar R-\bar\rho)^n}\right) 
\end{eqnarray*}
for every $i\in \mathbb{N}$, {with a constant $\mathfrak{C}=\mathfrak{C}(n,\alpha,\lambda,\Lambda, {H_n})$}. {For brevity, from now on, we will still denote with the same symbol $\mathfrak{C}$ the different constants varying from line to line. Only relevant changes will be suitably highlighted.}
Therefore, we can let $i\to \infty$ in \eqref{stimapreitepassoin} thus getting
\begin{eqnarray}\label{stimaquasifin}
	&&||V(Du)||_{L^{\infty}(B_{\bar\rho})}\cr\cr
    &\le&  \mathfrak{C}\left(\Big(1+\log^2\|V(Du)\|_{\infty}\Big)^\frac{n}{2}\exp\Big( \frac{n}{2}  (1 + \log^2 \|{V(Du)}\|_{\infty})^{\frac{1}{\alpha}}\Big)+\frac{\Big(1+\log^2\|V(Du)\|_{\infty}\Big)^\frac{n}{2}}{(\bar R-\bar\rho)^n}\right)   \cr\cr 
 &&\qquad\cdot\int_{B_{ \bar R}}V(Du)
\, \dx.
\end{eqnarray}
If $\|V(Du)\|_{\infty}\le e$, the proof ends here.
\\ Hence from now on, we shall suppose that $\|V(Du)\|_{\infty}> e$.
 {In this case, using \eqref{ordinelognuovo} to estimate the first term in the right hand side of \eqref{stimaquasifin} which is legitimate since $\alpha>4$, we get
\begin{eqnarray}\label{stimaquasifin1}
	||V(Du)||_{L^{\infty}(B_{\bar\rho})} &\le& 
 \mathfrak{C}\left( \frac{ \|V(Du)\|_{\infty} }{\log\big(\|V(Du)\|_{\infty}\big)}+\frac{\Big(1+\log^2\|V(Du)\|_{\infty}\Big)^\frac{n}{2}}{( \bar R- \bar\rho )^{n}}+\hat{\mathfrak{C}}\right)\cr\cr
&&\qquad\cdot\int_{B_{ \bar R}}V(Du)\, \dx.
    \end{eqnarray}
Using Lemma \ref{PdN} 
 with $$\sigma=\frac{1}{4\log(1+e)},\quad \alpha=\beta=1,\quad  s=\frac{ \|V(Du)\|_{\infty} }{\log\big(\|V(Du)\|_{\infty}\big)}\,\, \text{and}\,\, t=\mathfrak{C}\int_{B_{ \bar R}}V(Du)
\, \dx$$
  yields
  \begin{eqnarray}
   &&\mathfrak{C}\frac{ \|V(Du)\|_{\infty} }{\log\big(\|V(Du)\|_{\infty}\big)} \int_{B_{ \bar R}}V(Du)\,\dx\le \frac{1}{4\log(1+e)}\frac{ \|V(Du)\|_{\infty} }{\log\big(\|V(Du)\|_{\infty}\big)}\log\left(1+\frac{ \|V(Du)\|_{\infty} }{\log\big(\|V(Du)\|_{\infty}\big)}\right)\cr\cr
&&\qquad\qquad+\mathfrak{C}\int_{B_{ \bar R}}V(Du)
\, \dx\cdot\exp\left(\tilde{\mathfrak{C}}\int_{B_{ \bar R}}V(Du)
\, \dx\right)\cr\cr
&\stackrel{\eqref{elem}}{\le}& \frac{1}{4}\|V(Du)\|_{\infty}+\mathfrak{C}\int_{B_{ \bar R}}V(Du)
\, \dx\cdot\exp\left(\tilde{\mathfrak{C}}\int_{B_{ \bar R}}V(Du)
\, \dx\right).
  \end{eqnarray}
  Inserting this inequality in \eqref{stimaquasifin1}, we deduce that
\begin{eqnarray*}
	||V(Du)||_{L^{\infty}(B_{\bar\rho})} &\le&  \frac14||V(Du)||_{L^{\infty}(B_{ \bar R})}+ \mathfrak{C}\int_{B_{ \bar R}}V(Du)
\, \dx\cdot\exp\Bigg(\tilde{\mathfrak{C}}\int_{B_{ \bar R}}V(Du)\, \dx\Bigg)\cr\cr
&&+\mathfrak{C}\frac{\log^n\|V(Du)\|_{\infty}}{( \bar R- \bar\rho )^{n}}\int_{B_{ \bar R}}V(Du)\, \dx+\tilde{\mathfrak{C}}\left(1+\frac{1}{( \bar R- \bar\rho )^{n}}\right)\int_{B_{ \bar R}}V(Du)\, \dx.
\end{eqnarray*}
Using  \eqref{elembis} 
in the third term of the right hand side of the previous inequality with $$t=\log^n\|V(Du)\|_{\infty}\qquad \text{and} \qquad s=\frac{\mathfrak{C}}{( \bar R- \bar\rho )^{n}}\int_{B_{ \bar R}}V(Du)\, \dx$$
we obtain
\begin{eqnarray*}
	||V(Du)||_{L^{\infty}(B_{\bar\rho})} &\le&  \frac12||V(Du)||_{L^{\infty}(B_{ \bar R})}+ \mathfrak{C}\int_{B_{ \bar R}}V(Du)
\, \dx\cdot\exp\Bigg(\tilde{\mathfrak{C}}\int_{B_{ \bar R}}V(Du)\, \dx\Bigg)\cr\cr
    &&+\left(\frac{\mathfrak{C}}{( \bar R- \bar\rho )^{n}}\int_{B_{ \bar R}}V(Du)\, \dx\right)\log^n\left(1+\frac{\mathfrak{C}}{( \bar R- \bar\rho )^{n}}\int_{B_{ \bar R}}V(Du)\, \dx\right) \cr\cr
 &&+\frac{\mathfrak{C}}{( \bar R- \bar\rho )^{n}}\int_{B_{ \bar R}}V(Du)\, \dx+\hat{\mathfrak{C}}\int_{B_{ \bar R}}V(Du)\, \dx,
\end{eqnarray*}
that obviously implies
\begin{eqnarray*}
	||V(Du)||_{L^{\infty}(B_{\bar\rho})} &\le&  \frac12||V(Du)||_{L^{\infty}(B_{ \bar R})}+\mathfrak{C}\int_{B_{ \bar R}}V(Du)
\, \dx\cdot \exp\Bigg(\tilde{\mathfrak{C}}\int_{B_{ R_0}}V(Du)\, \dx\Bigg)\cr\cr
    &&+\left(\frac{\mathfrak{C}}{( \bar R- \bar\rho )^{n}}\int_{B_{ R_0}}V(Du)\, \dx\right)\log^n\left(1+\frac{\mathfrak{C}}{( \bar R- \bar\rho )^{n}}\int_{B_{  R_0}}V(Du)\, \dx\right) \cr\cr
 &&+\frac{\mathfrak{C}}{( \bar R- \bar\rho )^{n}}\int_{B_{ R_0}}V(Du)\, \dx+\hat{\mathfrak{C}}\int_{B_{ R_0}}V(Du)\, \dx.
\end{eqnarray*}
Since previous estimate holds true for every $\rho_0<\bar\rho<\bar R<R_0$, the iteration Lemma \ref{lem:Giaq} implies the conclusion.}
 \end{proof}

{\subsection{Proof of Theorem \ref{thmapriori} in the case $n = 2$}}
 
 \begin{proof} Let us fix a ball $B_{R_0}\Subset \Omega$, radii $\rho_{0}<\bar\rho<\rho<R<\bar R<R_{0}$
 and a cut-off function $\eta\in C^1_0(B_{R})$ such that $\eta=1$ on $B_{\rho}$ and $|D\eta|\le \frac{C}{R-\rho}$. The main change with respect to the proof of the case $n\ge 3$ is the use of the Sobolev inequality at Theorem \ref{sobtru} with the function $w=\eta(1+|Du|^2)^{\frac{\gamma + 1}{4}} = \eta V(Du)^{\frac{\gamma + 1}{2}}$ and $n=2$, where $V(Du)$ has been defined at  \eqref{V(Du)}, i.e.
 \begin{eqnarray*}
&&{\Big|\Big|\eta V(Du)^{\frac{\gamma + 1}{2}}\Big|\Big|^2_{{\rm EXP}_2(B_R)}}\le \mathcal{C}_S\int_{B_R}\left|D\left(\eta V(Du)^{\frac{\gamma + 1}{2}} \right)\right|^2\, \dx\cr\cr
&\le& C \, (1+\gamma)^2\int_{B_R}\eta^2 V(Du)^{\gamma - 1}|D^2 u|^2+\frac{C}{(R-\rho)^2}\int_{B_R} V(Du)^{\gamma + 1}\, \dx.
 \end{eqnarray*}
 Using \eqref{estprimadellasob} to estimate the first term in the right hand side of previous inequality we get
 \begin{eqnarray*}
&&{\Big|\Big|\eta V(Du)^{\frac{\gamma+1}{2}}\Big|\Big|^2_{{\rm EXP}_2(B_R)}}\cr\cr
&\le& C \, (1+\gamma)^4   \, \int_{B_{R}}\eta^2  h^2 V(Du)^{1 + \gamma} (1 + \log^2(V(Du)) ) \, \dx\cr\cr
&&\quad+  \frac{C \, (1+\gamma)^2}{(R-\rho)^2}  \, \int_{B_{R}}   V(Du)^{1 + \gamma} (1+\log V(Du))  \, \dx.
\end{eqnarray*} 
 The a priori assumption $Du\in L^{\infty}_{\mathrm{loc}}(\Omega)$ implies
\begin{eqnarray}\label{primastiman=2}
&&\Big|\Big|\eta V(Du)^{\frac{1+\gamma}{2}}\Big|\Big|^2_{{\rm EXP}_2(B_R)}\cr\cr
&\le& C \, (1+\gamma)^4 \big(1 + \log^2||V(Du)||_{\infty}\big) \, \int_{B_{R}}\eta^2  h^2 V(Du)^{1 + \gamma}  \, \dx\cr\cr
&&\quad + \frac{C(1+\gamma)^2\big(1 + \log||V(Du)||_{\infty}\big)}{(R-\rho)^2}  \, \int_{B_{R}}   V(Du)^{1 + \gamma}  \, \dx\cr\cr
&\le& C \, (1+\gamma)^4 \big(1 + \log^2||V(Du)||_{\infty}\big) \, \int_{B_{R}\cap\{h\le \varepsilon\}}\eta^2  h^2 V(Du)^{1 + \gamma}  \, \dx\cr\cr
&&+ C \, (1+\gamma)^4 \big(1 + \log^2||V(Du)||_{\infty}\big)  \, \int_{B_{R}\cap\{h> \varepsilon\}}\eta^2  h^2 V(Du)^{1 + \gamma}  \, \dx\cr\cr
&&\quad+  \frac{C \, (1+\gamma)^2\big(1 + \log||V(Du)||_{\infty}\big)}{(R-\rho)^2}  \, \int_{B_{R}}   V(Du)^{1 + \gamma}  \, \dx\cr\cr
&\le& C \,\varepsilon^2(1+\gamma)^4 \big(1 + \log^2||V(Du)||_{\infty}\big) \int_{B_{R}}\eta^2   V(Du)^{1 + \gamma}  \, \dx\cr\cr
&&+ C \, (1+\gamma)^4 \big(1 + \log^2||V(Du)||_{\infty}\big) \, \int_{B_{R}}\eta^2  h^2\chi_{\{h> \varepsilon\}} V(Du)^{1 + \gamma}  \, \dx\cr\cr
&&\quad+  \frac{C \, (1+\gamma)^2\big(1 + \log||V(Du)||_{\infty}\big)}{(R-\rho)^2}  \, \int_{B_{R}}   V(Du)^{1 + \gamma}  \, \dx,
\end{eqnarray}
where $\varepsilon>0$ will be chosen later and we shortened with $||V(Du)||_{\infty}$ the norm $||V(Du)||_{L^\infty(B_{\bar R})}$.
Using Proposition \ref{HoldOrlicz} to estimate the second integral in the right hand side of the previous estimate, we obtain
\begin{eqnarray}\label{primastiman=2}
&&{\Big|\Big|\eta V(Du)^{\frac{1+\gamma}{2}}\Big|\Big|^2_{{\rm EXP}_2(B_R)}}\cr\cr
&\le&  C\,\varepsilon^2(1+\gamma)^4 \big(1 + \log^2||V(Du)||_{\infty}\big) \int_{B_{R}}\eta^2   V(Du)^{1 + \gamma}  \, \dx\cr\cr
&&+ C (1+\gamma)^4 \big(1 + \log^2\||V(Du)||_{\infty}\big)\big|\big|h^2\chi_{\{h> \varepsilon\}} \big|\big|_{L\log L(B_R)}{\Big|\Big|\eta^2 V(Du)^{{1+\gamma}}\Big|\Big|_{{\rm EXP}(B_R)}}\cr\cr
&&\quad+  \frac{C (1+\gamma)^2\big(1 + \log||V(Du)||_{\infty}\big)}{(R-\rho)^2}  \, \int_{B_{R}}   V(Du)^{1 + \gamma}  \, \dx\cr\cr
&\le& C\,\varepsilon^2(1+\gamma)^4\big(1 + \log^2||V(Du)||_{\infty}\big)\int_{B_{R}}\eta^2   V(Du)^{1 + \gamma}  \, \dx\cr\cr
&&\stackrel{\eqref{galquadrato}}{+} C(1+\gamma)^4 \big(1 + \log^2||V(Du)||_{\infty}\big)\big|\big|h^2\chi_{\{h> \varepsilon\}} \big|\big|_{L\log L({B_{R_0}})}{\Big|\Big|\eta V(Du)^{^{\frac{1+\gamma}{2}}}\Big|\Big|^2_{{\rm EXP}_{2}(B_R)}}\cr\cr
&&\quad+  \frac{C(1+\gamma)^2\big(1 + \log||V(Du)||_{\infty}\big)}{(R-\rho)^2}  \, \int_{B_{R}}   V(Du)^{1 + \gamma}  \, \dx.
\end{eqnarray}
 Now, using \eqref{alphallogl} {and the fact that 
 \(
\log(e + h^2) \le \, \log ((e + h)^2) = 2 \log(e + h), 
 \)} 
 we get
\begin{eqnarray*}
 &&\big|\big|h^2\chi_{\{h> \varepsilon\}}\big|\big|_{L\log L({B_{R_0}})}\cr\cr
 &\le& c({\beta,R_0})\left(\int_{{B_{R_0}}}h^2\chi_{\{h> \varepsilon\}}\log (e+h^2\chi_{\{h>\varepsilon\}})\,\dx\right)^{\frac{1}{{\beta}}} + {\int_{{B_{R_0}}}h^2\chi_{\{h> \varepsilon\}}\log (e+h^2\chi_{\{h> \varepsilon\}})\,\dx}\cr\cr  &=& c({\beta},R_0) \left(\int_{B_{{B_{R_0}}}\cap\{h> \varepsilon\}}h^2\log (e+h^2)\,\dx\right)^{\frac{1}{{\beta}}} + {\int_{B_{{B_{R_0}}}\cap\{h> \varepsilon\}}h^2\log (e+h^2)\,\dx}\cr\cr
 &\le& c({\beta,R_0})\left(\int_{B_{{B_{R_0}}}\cap\{h> \varepsilon\}}h^2\frac{\log^{2\alpha} (e+h)}{\log^{2\alpha-1}(e+h)}\,\dx\right)^{\frac{1}{{\beta}}} {+ \int_{B_{{B_{R_0}}}\cap\{h> \varepsilon\}}h^2\frac{\log^{2\alpha} (e+h)}{\log^{2\alpha-1}(e+h)}\,\dx} \cr\cr
 &\le& \frac{c({\beta,R_0})}{\log^{\frac{2\alpha-1}{\beta}}(e+\varepsilon)}\left(\int_{B_{{B_{R_0}}}\cap\{h> \varepsilon\}}h^2\log^{2\alpha}(e+h)\,\dx\right)^{\frac{1}{{\beta}}}\cr\cr
 && + {\frac{1}{\log^{2\alpha-1}(e+\varepsilon)}\left(\int_{B_{{B_{R_0}}}\cap\{h> \varepsilon\}}h^2\log^{2\alpha}(e+h)\,\dx\right),}
\end{eqnarray*}
which is finite by virtue of assumption {(F4)} for $\alpha>\frac92$ {(this bound on $\alpha$ will be exploited later)}. {Since $\log(e + \varepsilon) \ge \, 1$ and $\beta > 1,$ we have that
\[
\log^{\frac{2 \alpha - 1}{\beta}}(e + \varepsilon)\le \, \log^{2 \alpha - 1}(e + \varepsilon).
\]
Hence,} using the previous estimate in \eqref{primastiman=2}, we obtain
\begin{eqnarray}\label{primastiman=2bis}
&&{\Big|\Big|\eta V(Du)^{\frac{1+\gamma}{2}}\Big|\Big|^2_{{\rm EXP}_2(B_R)}}\cr\cr
&\le&  C\varepsilon^2(1+\gamma)^4 \big(1+\log^2||V(Du)||_{\infty}\big)   \int_{B_{R}}\eta^2   V(Du)^{1 + \gamma}  \, \dx\cr\cr
&&+\frac{{C(\beta,R_0)}(1+\gamma)^4}{\log^{\frac{2\alpha-1}{{\beta}}}(e+\varepsilon)} \big(1+\log^2||V(Du)||_{\infty}\big) \cr\cr
&& \cdot\left [\left(\int_{B_R}h^2\log^{2\alpha}(e+h)\, \dx\right)^{\frac{1}{{\beta}}}
+ {\int_{B_R}h^2\log^{2\alpha}(e+h)\, \dx}
\right ]
\Big|\Big|\eta V(Du)^{^{\frac{1+\gamma}{2}}}\Big|\Big|^2_{{\rm EXP}_2(B_R)}\cr\cr 
&&\quad+  \frac{{C(1+\gamma)^2}\big(1+\log||V(Du)||_{\infty}\big)}{(R-\rho)^2}  \, \int_{B_{R}}   V(Du)^{1 + \gamma}  \, \dx\cr\cr
&\le&  C\varepsilon^2(1+\gamma)^4 \big(1+\log^2||V(Du)||_{\infty}\big)   \int_{B_{R}}\eta^2   V(Du)^{1 + \gamma}  \, \dx\cr\cr
&&+\frac{{C(\beta,R_0)[H_2^{\frac{2}{\beta}} + H_2^2]}(1+\gamma)^4}{\log^{\frac{2\alpha-1}{{\beta}}}(e+\varepsilon)} \big(1+\log^2||V(Du)||_{\infty}\big)\Big|\Big|\eta V(Du)^{^{\frac{1+\gamma}{2}}}\Big|\Big|^2_{{\rm EXP}_2(B_R)}\cr\cr 
&&\quad + \frac{C(1+\gamma)^2\big(1+\log||V(Du)||_{\infty}\big)}{(R-\rho)^2}  \, \int_{B_{R}}   V(Du)^{1 + \gamma}  \, \dx,
\end{eqnarray}
where, we used as before the notation
\[
{H_2} := \left (\int_{B_{R_0}} h^2 \log^{2\alpha}(e+h) \, \dx\right )^{{\frac{1}{2}}}.
\]
Choosing now $\varepsilon>0$ such that
$$\frac{{C(\beta,R_0)}{[H_2^{\frac{2}{\beta}} + H_2^2]}(1+\gamma)^4}{\log^{\frac{2\alpha-1}{{\beta}}}(e+\varepsilon)} \big(1+\log^2||V(Du)||_{\infty}\big)=\frac{1}{2}$$  $$ \,\,\Longleftrightarrow\,\, e+\varepsilon=\exp\Big(C^{\frac{{\beta}}{2\alpha-1}}\,2^{\frac{{\beta}}{2\alpha-1}}{[H_2^{\frac{2}{\beta}} + H_2^2]^{\frac{\beta}{2\alpha-1}}}(1+\gamma)^{{\frac{4{\beta}}{2\alpha-1}}} \big(1+\log^2||V(Du)||_{\infty}\big)^{\frac{{\beta}}{2\alpha-1}}\Big)$$
estimate \eqref{primastiman=2bis} becomes
\begin{eqnarray}\label{primastiman=2ter}
&&{\Big|\Big|\eta V(Du)^{\frac{1+\gamma}{2}}\Big|\Big|^2_{{\rm EXP}_2(B_R)}}\le \frac{1}{2}{\Big|\Big|\eta V(Du)^{\frac{1+\gamma}{2}}\Big|\Big|^2_{{\rm EXP}_2(B_R)}}\cr\cr
&& + C\Big[(1+\gamma)^4 \big(1+\log^2||V(Du)||_{\infty}\big)\exp\Big(\mathfrak{C}(1+\gamma)^{{\frac{4{\beta}}{2\alpha-1}}} \big(1+\log^2||V(Du)||_{\infty}\big)^{\frac{{\beta}}{2\alpha-1}}\Big)\Big] \cr\cr &&\qquad\qquad \cdot \int_{B_{R}}\eta^2   V(Du)^{1 + \gamma}  \, \dx\cr\cr
&&\quad+  \frac{C(1+\gamma)^2\big(1+\log||V(Du)||_{\infty}\big)}{(R-\rho)^2} \, \int_{B_{R}}   V(Du)^{{1 + \gamma}}  \, \dx,
\end{eqnarray}
{where, for the sake of simplicity, we set
\[
\mathfrak{C} := C^{\frac{{\beta}}{2\alpha-1}}\,2^{\frac{{\beta}}{2\alpha-1}+1}{[H_2^{\frac{2}{\beta}} + H_2^2]^{\frac{\beta}{2\alpha-1}}}.
\]}
Reabsorbing the first integral in the right hand side by the left hand side, we infer
\begin{eqnarray}
&&  \frac{1}{2}{\Big|\Big|\eta V(Du)^{\frac{1+\gamma}{2}}\Big|\Big|^2_{{\rm EXP}_2(B_R)}}\cr\cr 
&\le & C\Big[(1+\gamma)^4 \big(1+\log^2||V(Du)||_{\infty}\big)\exp\Big(\mathfrak{C}(1+\gamma)^{{\frac{4{\beta}}{2\alpha-1}}} \big(1+\log^2||V(Du)||_{\infty}\big)^{\frac{{\beta}}{2\alpha-1}}\Big)\Big] \cr\cr&&\qquad\qquad\cdot  \int_{B_{R}}\eta^2   V(Du)^{1 + \gamma}  \, \dx\cr\cr
&&\quad  +\frac{{C(1+\gamma)^2\big(1+\log||V(Du)||_{\infty}\big)}}{(R-\rho)^2}  \, \int_{B_{R}}   V(Du)^{1 + \gamma}  \, \dx.
\end{eqnarray}
By Proposition \ref{normalpexp}, 
 recalling that $\eta=1$ on $B_\rho$, we get
 \begin{eqnarray*}\label{terzastiman=2}
&&\left(\int_{B_{\rho}}V(|Du|)^{(1+\gamma)\frac{2^*}{2}} \,\dx\right)^{\frac{2}{\cdot2^*}}\cr\cr
&\le & C(1+\gamma)^4 (1+\log^{2}||V(Du)||_{\infty}\big)\exp\Big(\mathfrak{C}(1+\gamma)^{{\frac{4{\beta}}{2\alpha-1}}} \big(1+\log^2||V(Du)||_{\infty}\big)^{\frac{{\beta}}{2\alpha-1}}\Big)\cr\cr     &&\qquad\quad \cdot\int_{B_{R}}\eta^2   V(Du)^{1 + \gamma}  \, \dx  +\frac{C(1+\gamma)^2(1+\log ||V(Du)||_{\infty})}{(R-\rho)^2}  \, \int_{B_{R}}   V(Du)^{1 + \gamma}  \, \dx\cr\cr 
\cr\cr
&\le & C(1+\gamma)^4\exp\Big(\mathfrak{C}(1+\gamma)^{{\frac{8{\beta}}{2\alpha-1}}} \Big) (1+\log^{2}||V(Du)||_{\infty}\big)\exp\Big( \big(1+\log^2||V(Du)||_{\infty}\big)^{\frac{{2\beta}}{2\alpha-1}}\Big)\cr\cr     &&\qquad\quad \cdot\int_{B_{R}}\eta^2   V(Du)^{1 + \gamma}  \, \dx  +\frac{C(1+\gamma)^2(1+\log ||V(Du)||_{\infty})}{(R-\rho)^2}  \, \int_{B_{R}}   V(Du)^{1 + \gamma}  \, \dx,
 \end{eqnarray*}
where we used the elementary inequality $\exp(st)\le \exp(s^2)\exp(t^2)$ and  now $2^*$ denotes any finite exponent greater than $2$, since $n=2$. 
As in the case $n>2$,  considering  the decreasing sequence of  radii defined at \eqref{raggi}
and the  increasing sequence of exponents defined at \eqref{esponenti},
 we may iterate estimate \eqref{terzastiman=2}  on the concentric balls $B_{\rho_i}$  with exponents $p_i$,
thus obtaining
\begin{eqnarray}\label{stimapreitepassoi}
&&\left ( \int_{B_{\rho_{i+1}}}  V(Du)^{p_{i+1}}\, \dx\right )^{\frac{1}{p_{i+1}}} 	\cr\cr
&\le&\!\! \displaystyle{\prod_{j=0}^{i}}\Bigg(Cp_j^{4}\exp\Big(\mathfrak{C}p_j^{{\frac{8{\beta}}{2\alpha-1}}}\Big)\big(1+\log^{2}||V(Du)||_{\infty}\big)\exp\Big(\big(1+\log^2||V(Du)||_{\infty}\big)^{\frac{{2\beta}}{2\alpha-1}}\Big)\cr\cr  &&\qquad\qquad+C\frac{p_j^2(1+\log ||V(Du)||_{\infty})}{(\rho_j - \rho_{j+1})^{2}}  \Bigg)^{\frac{1}{p_j}}\cdot\int_{B_{\bar R}}V(Du) \, \dx\cr\cr
&=& \displaystyle{\prod_{j=0}^{i}} \Bigg(Cp_j^4\exp\Big(\mathfrak{C}p_j^{{\frac{8{\beta}}{2\alpha-1}}}\Big)\big(1+\log^{2}||V(Du)||_{\infty}\big)\exp\Big(\big(1+\log^2||V(Du)||_{\infty}\big)^{\frac{{2\beta}}{2\alpha-1}}\Big)\cr\cr  &&\qquad\qquad+C\frac{4^{j+1}p_j^2(1+\log ||V(Du)||_{\infty})}{(\bar R - \bar\rho)^2}  \Bigg)^{\frac{1}{p_j}}\cdot\int_{B_{\bar R}}V(Du)  \, \dx\cr\cr
&\le & \displaystyle{\prod_{j=0}^{i}} \left((1+\log^2||V(Du)||_{\infty})\exp\Big( \big(1+\log^2||V(Du)||_{\infty}\big)^{\frac{{2\beta}}{2\alpha-1}}\Big) +\frac{(1+\log ||V(Du)||_{\infty})}{(\bar R - \bar\rho)^2}   \right)^{\frac{1}{p_j}}\cr\cr
&&\qquad\qquad\left(Cp_j^4\exp\Big(\mathfrak{C}p_j^{{\frac{8{\beta}}{2\alpha-1}}}\Big) +C4^{j+1}p_j^2 \right)^{\frac{1}{p_j}}\cdot\int_{B_{\bar R}}V(Du)  \, \dx\cr\cr
&&=:\tilde\Theta_i \int_{B_{\bar R}}V(Du)\, \dx .
\end{eqnarray}
Since
$$\tilde\Theta_i \le \exp\Bigg(\sum_{j=0}^i \frac{1}{p_j}\log\Bigg(\big(1+\log^{2}||V(Du)||_{\infty}\big)\exp\Big(\big(1+\log^2||V(Du)||_{\infty}\big)^{\frac{{2\beta}}{2\alpha-1}}\Big)$$$$\qquad\qquad  +\frac{(1+\log ||V(Du)||_{\infty})}{(\bar R - \bar\rho)^2}\Bigg) \Bigg) \cdot \exp\left(\sum_{j=0}^i \frac{1}{p_j}\log\left(C p_j^{4}\exp\Big(\mathfrak{C}p_j^{{\frac{8{\beta}}{2\alpha-1}}}\Big)+C4^{j+1}p_j^2 \right)\right)$$
$$\le \exp\Bigg(\sum_{j=0}^\infty \frac{1}{p_j}\log\Bigg(\big(1+\log^{2}||V(Du)||_{\infty}\big)\exp\Big(\big(1+\log^2||V(Du)||_{\infty}\big)^{\frac{{2\beta}}{2\alpha-1}}\Big)$$$$\qquad\qquad  +\frac{(1+\log ||V(Du)||_{\infty})}{(\bar R - \bar\rho)^2}\Bigg) \Bigg)  \cdot \exp\left(\sum_{j=0}^\infty \frac{1}{p_j}\log\left(Cp_j^4\exp\Big(\mathfrak{C}p_j^{{\frac{8{\beta}}{2\alpha-1}}}\Big) +C4^{j+1}p_j^2 \right)\right)$$
$$\le \Bigg(\big(1+\log^{2}||V(Du)||_{\infty}\big)\exp\Big(\big(1+\log^2||V(Du)||_{\infty}\big)^{\frac{{2\beta}}{2\alpha-1}}\Big) +\frac{(1+\log ||V(Du)||_{\infty})}{(\bar R - \bar\rho)^2}\Bigg)^{\frac{ 2^*}{2^*-2}}$$
$$ \cdot \exp\left(\sum_{j=0}^\infty \frac{1}{p_j}\left(\log\left({2}Cp_j^4 \right)+\log\left(2\exp\Big(\mathfrak{C}p_j^{{\frac{8{\beta}}{2\alpha-1}}}\Big)  \right)+\log(2C4^{j+1})+\log\left(2p_j^2\right)\right)\right)$$
$$=\Bigg(\big(1+\log^{2}||V(Du)||_{\infty}\big)\exp\Big(\big(1+\log^2||V(Du)||_{\infty}\big)^{\frac{{2\beta}}{2\alpha-1}}\Big) +\frac{(1+\log ||V(Du)||_{\infty})}{(\bar R - \bar\rho)^2} \Bigg)^{\frac{ 2^*}{2^*-2}}$$ $$\cdot\tilde\Theta_1\cdot \tilde\Theta_2\cdot \tilde\Theta_3\cdot \tilde\Theta_4,$$
where we used that, by the definition of $p_j$, it holds $\sum_{j=0}^{{\infty}}\frac{1}{p_j}= \frac{2^*}{ 2^*-2} $. 
Since $\alpha>\frac{9}{2}$, we may choose $1<{\beta}<\frac{2\alpha-1}{8}$, thus getting
$$ \tilde\Theta_2 \le C(n,\alpha,\mathfrak{C})$$
and one can easily check that
$$ \tilde\Theta_1\cdot \tilde\Theta_3\cdot \tilde\Theta_4 \le C . $$
Combining  the estimates above, we conclude that
\begin{eqnarray*}
\tilde \Theta_i&\le& \tilde{\mathfrak{C}}\Bigg(\big(1+\log^{2}||V(Du)||_{\infty}\big)\exp\Big(\big(1+\log^2||V(Du)||_{\infty}\big)^{\frac{{2\beta}}{2\alpha-1}}\Big) +\frac{(1+\log ||V(Du)||_{\infty})}{(\bar R - \bar\rho)^2} \Bigg)^{\frac{ 2^*}{2^*-2}}
\end{eqnarray*}
for every $i\in \mathbb{N}$ and with $ \tilde{\mathfrak{C}}= \tilde{\mathfrak{C}}(n,\alpha, {H_2},\lambda,\Lambda)$.
Therefore, we can let $i\to \infty$ in \eqref{stimapreitepassoi} thus getting
\begin{eqnarray}
    ||V(Du)||_{L^\infty(B_{\bar\rho})}&\le& \mathfrak{C}\Bigg(\big(1+\log^{2}||V(Du)||_{\infty}\big)\exp\Big(\big(1+\log^2||V(Du)||_{\infty}\big)^{\frac{{2\beta}}{2\alpha-1}}\Big)\cr\cr
    &&\qquad+\frac{(1+\log^{{2}} ||V(Du)||_{\infty})}{(\bar R - \bar\rho)^2} \Bigg)^{\frac{ 2^*}{2^*-2}}\cdot\int_{B_{\bar R}}V(Du)  \, \dx\cr\cr 
    &\le& \mathfrak{C}\Bigg(\big(1+\log^2||V(Du)||_{\infty}\big)^{\frac{  2^*}{2^*-2}}\exp\Big({\frac{ 2^*}{2^*-2}}\big(1+\log^2||V(Du)||_{\infty}\big)^{\frac{{2\beta}}{2\alpha-1}}\Big)\cr\cr
    &&\qquad+\frac{(1+\log^{\frac{{2^* 2}}{2^*-2}} ||V(Du)||_{\infty})}{(\bar R - \bar\rho)^{\frac{ 2\cdot 2^*}{2^*-2}}} \Bigg)\cdot\int_{B_{\bar R}}V(Du)  \, \dx
\end{eqnarray}
Since by our choice of ${\beta}$, it holds $\frac{2\alpha-1}{2{\beta}}>2$, we may use \eqref{ordinelognuovo} to deduce
\begin{eqnarray}
    ||V(Du)||_{L^\infty(B_{\bar\rho})}&\le& C\left(\frac{||V(Du)||_{L^\infty(B_{\bar R})}}{\log \big(||V(Du)||_{L^\infty(B_{\bar R})}\big)}+\frac{(1+\log^{\frac{ {2^* 2}}{2^*-2}} ||V(Du)||_{\infty})}{(\bar R - \bar\rho)^{\frac{ 2\cdot 2^*}{2^*-2}}}+1\right)\cr\cr 
&&\qquad\qquad\cdot\int_{B_{\bar R}}V(Du)  \, \dx.
\end{eqnarray}
At this point the proof can be concluded as in the case {$n\ge3$ with the obvious modifications}.
\end{proof}

 \section{The approximation}

\label{quattro}

In this section we give only a sketch of  the approximation procedure that comes along quite standard arguments.
\\
For a smooth mollifier $\varphi\in C^\infty_0(B_1(0))$ and $\varepsilon\in (0,1]$, let us consider 
$$F_{\varepsilon}(x,\xi)=\int_{B_{1}(0)} F(x+\varepsilon y, \xi)\varphi(y)\,dy.$$
It is easy to check that $F_{\varepsilon}( x,\xi)$ still satisfies 
\[
{c} |\xi|\log(1 + |\xi|) \le \, F_{\varepsilon}(x, \xi)  \le \, {C} |\xi|\log(1 + |\xi|) \tag{F1'}
\]
\[
\lambda\frac{|\mu|^2}{1+|\xi|}
\le
\langle {D_{{\xi\xi}}F_{\varepsilon}}(x,\xi)\mu,\mu\rangle
\le
\Lambda
\frac{\log(1+|\xi|)}{|\xi|}
|\mu|^2,
\tag{F2'}
\]  
\[\vert D_{\xi x} F_{\varepsilon}(x,\xi)\vert \le h_\varepsilon(x) \log(1+|\xi |) \tag{F3'},
\]
with constants independent of $\varepsilon,$ where $h_\varepsilon$ is the usual mollification of $h$ through $\varphi$ and it also satisfies assumption (F4) by
Jensen’s inequality.
\color{black}
Let $u$ be a local minimizer  of the functional {$\mathcal{F}$} and let $v_{\varepsilon}\in u+W^{1,L\log L}_0(B_{R_0})$ be the unique solution  to the problem
$$\min\left\{\mathcal{F}_{\varepsilon} (v,B_{R_0}):\,\,\, v\in u+W^{1,L\log L}_0(B_{R_0})\right\},$$
where $$\mathcal{F}_{\varepsilon} (v,B_{R_0})=\int_{B_{R_0}}F_{\varepsilon}(x,Dv)\,\dx.$$
Hence, we  may use the regularity results in see  \cite{Giova, BR, G20, BDPdN} to deduce that  $v_{\varepsilon}\in W^{1,\infty}_{\rm loc}(B_{R_0})\cap W^{2,2}_{\rm loc}(B_{R_0}) $   for every $\varepsilon>0$. Therefore, we are legitimate to use the a priori estimate of Theorem \ref{thmapriori} to deduce that
\begin{eqnarray}\label{aprioriapp}
	||V(Dv_{\varepsilon})||_{L^{\infty}(B_{\rho})} &\le& 
  \mathfrak{C}\int_{B_{  R}}V(Dv_{\varepsilon})
\, \dx\cdot \exp\Bigg(\tilde{\mathfrak{C}}\int_{B_{ R}}V(Dv_{\varepsilon})\, \dx\Bigg)\cr\cr
    &&+\left(\frac{\mathfrak{C}}{(  R- \rho )^{{\theta}}}\int_{B_{ R}}V(Dv_{\varepsilon})\, \dx\right)\log^{{\theta}}\left(1+\frac{\mathfrak{C}}{(  R- \rho )^{{\theta}}}\int_{B_{  R}}V(Dv_{\varepsilon})\, \dx\right) \cr\cr
 &&+\frac{\mathfrak{C}}{(  R- \rho )^{{\theta}}}\int_{B_{ R}}V(Dv_{\varepsilon})\, \dx+\hat{\mathfrak{C}}\int_{B_{ R}}V(Dv_{\varepsilon})\, \dx,
\end{eqnarray}
for every $B_\rho\subset B_R\Subset B_{R_0}$. On the other hand, the left inequality in {(F1')} implies that
\begin{eqnarray}\label{stima1}
&& {c} \int_{B_{R_0}}|Dv_{\varepsilon}|\log(1+|Dv_{\varepsilon}|)\, \dx\le \int_{B_{R_0}}F_{\varepsilon} (x,Dv_{\varepsilon})\,\dx \label{normaLp} \\[2mm]
&\le& \int_{B_{R_0}}F_{\varepsilon} (x,Du)\,\dx\le {C}\int_{B_{R_0}}|Du|\log (1+|Du|)\,\dx  \nonumber
\end{eqnarray}
where we used the minimality of $v_{\varepsilon}$ and the right inequality in {(F1').}
Therefore $\big(Dv_{\varepsilon}\big)_\varepsilon$ is an equi-integrable sequence in $u+W^{1,1}(B_{R_0})$ and so, by Dunford-Pettis  Theorem, there exists $v\in u+W^{1,1}(B_{R_0})$ such that \begin{equation}\label{weakLp}
  Dv_{\varepsilon} \rightharpoonup Dv \qquad \text{in}\,\, L^1(B_{R_0}).  
\end{equation}
The lower semicontinuity of functional  at \eqref{functional} implies
\begin{eqnarray}
&&\int_{B_{R_0}}|Dv|\log(1+|Dv|)\, \dx\le {c}\int_{B_{R_0}}F(x,Dv)\,\dx \cr\cr
&\le& {c} \liminf_{\varepsilon\to 0}\int_{B_{R_0}}F (x,Dv_{\varepsilon})\,\dx\le {\tilde{C}} \liminf_{\varepsilon\to 0}\int_{B_{R_0}}|Dv_\varepsilon|\log(1+|Dv_\varepsilon|)\, \dx\cr\cr &\le& {\tilde{C}} \int_{B_{R_0}}|Du|\log(1+|Du|)\,\dx, 
\end{eqnarray}
{being $\tilde{C} := c \, C$} and where in last inequality we used \eqref{stima1}.
Using {once more \eqref{stima1}} in \eqref{aprioriapp}, we get
\begin{eqnarray}\label{stimaconM}
||V(Dv_{\varepsilon})||_{L^\infty \left (B_{\frac{R_0}{4}}\right )}&\le&C\int_{B_{  R_0}}\Big(1+|Du|\log(1+|Du|)\Big)
\, \dx\cr\cr 
&&\qquad\cdot \exp\Bigg(C\int_{B_{ R_0}}\Big(1+|Du|\log(1+|Du|)\Big)\, \dx\Bigg)\cr\cr
    &&+\frac{C}{ R_0 ^{{\theta}}}\int_{B_{ R_0}}\Big(1+|Du|\log(1+|Du|)\Big)\, \dx\cr\cr    &&\qquad\cdot\log^{{\theta}}\left(1+\frac{C}{R_0^{{\theta}}}\int_{B_{  R}} \Big(1+|Du|\log(1+|Du|)\Big)\, \dx\right) \cr\cr
 &&+C\left(\frac{1}{ R_0^{{\theta}}}+1\right)\int_{B_{ R_0}}\Big(1+|Du|\log(1+|Du|)\Big)\, \dx \cr\cr
 &&{=: M.}
\end{eqnarray}
Therefore 
\begin{equation}\label{weakLpbi}
  v_{\varepsilon} \stackrel{*}\rightharpoonup v \qquad \text{weakly * in}\,\, W^{1,\infty}_{\mathrm{loc}}(B_{R_0})  
\end{equation}
and, by the lower semicontinuity of the norm, also
\begin{eqnarray}\label{apriori2bis}
||V(Dv)||_{L^\infty \left (B_{\frac{R_0}{4}}\right )}
&\le& C\int_{B_{  R_0}}\Big(1+|Du|\log(1+|Du|)\Big)
\, \dx\cr\cr 
&&\qquad\cdot \exp\Bigg(C\int_{B_{ R_0}}\Big(1+|Du|\log(1+|Du|)\Big)\, \dx\Bigg)\cr\cr
    &&+\frac{C}{ R_0 ^{{\theta}}}\int_{B_{ R_0}}\Big(1+|Du|\log(1+|Du|)\Big)\, \dx\cr\cr
&&\qquad\cdot\log^{{\theta}}\left(1+\frac{C}{R_0^{{\theta}}}\int_{B_{  R}} \Big(1+|Du|\log(1+|Du|)\Big)\, \dx\right) \cr\cr
 &&+C\left(\frac{1}{ R_0^{{\theta}}}+1\right)\int_{B_{ R_0}}|Du|\log(1+|Du|)\, \dx.
\end{eqnarray}
At this point, the conclusion of the proof, that is showing that $v=u$ a.e. in $B_{R_0}$, follows by standard arguments that can be found for instance in \cite{EMM16,CGGP,CGHPdN20}. 
\\
Indeed, by means of standard results about the convolutions and recalling \cite[Theorem 3]{Cianchi}, it holds that, for every $\delta > 0,$ for $\varepsilon$ sufficiently small 
\[
|F(x, \xi) - F_{\varepsilon}(x, \xi)| \le \, \delta  \qquad \textnormal{for a.a. $x \in \overline{B_{R_0}}$ and for every $|\xi| \in 1 + M \, \log (1 + M),$}
\]
being $M$ given by \eqref{stimaconM}; therefore, for any $r < R_0,$ by the minimality of $v_{\varepsilon}$ we have
\begin{eqnarray*}
&& \int_{B_{r}}F(x,Dv_{\varepsilon})\,\dx \cr \cr
&=& \int_{B_{r}} [F(x,Dv_{\varepsilon}) - F_{\varepsilon}(x,Dv_{\varepsilon}) ] \,\dx + \int_{B_{r}}  F_{\varepsilon}(x,Dv_{\varepsilon}) \,\dx\cr \cr
&\le& \delta \, |B_{R_0}| + \int_{B_{R_0}}  F_{\varepsilon}(x,Du) \,\dx.
\end{eqnarray*}
By the lower semicontinuity of $\mathcal{F}$, passing to the limit as $\varepsilon \rightarrow 0,$ we deduce
\begin{eqnarray*}
&& \int_{B_{r}}F(x,Dv)\,\dx \cr\cr
&\le& \liminf_{\varepsilon\to 0}\int_{B_{r}}F(x,Dv_{\varepsilon})\,\dx \cr \cr
&\le& \liminf_{\varepsilon\to 0}\int_{B_{r}} [F(x,Dv_{\varepsilon}) - F_{\varepsilon}(x,Dv_{\varepsilon}) ] \,\dx + \liminf_{\varepsilon\to 0} \int_{B_{r}}  F_{\varepsilon}(x,Dv_{\varepsilon}) \,\dx\cr \cr
&\le& \limsup_{\varepsilon\to 0}\int_{B_{R_0}} [F(x,Dv_{\varepsilon}) - F_{\varepsilon}(x,Dv_{\varepsilon}) ] \,\dx +  \limsup_{\varepsilon\to 0}\int_{B_{R_0}}  F_{\varepsilon}(x,Du) \,\dx\cr \cr
&=& \delta \, |B_{R_0}| + \int_{B_{R_0}}  F(x,Du) \,\dx.
\end{eqnarray*}
for every $\delta > 0$ and $r < R_0,$ and then, for $r \rightarrow R$ and $\delta \rightarrow 0$
\[
\int_{B_{R_0}}F(x,Dv)\,\dx \le \, \int_{B_{R_0}}  F(x,Du) \,\dx.
\]
Therefore $v$ and $u$ are two different minimizers of {$\mathcal{F}$} in $B_{R_0}$ and thus, by the strict convexity of the functional {$\mathcal{F}$} it is possible to conclude that $v=u$ a.e. in $B_{R_0}$. 
\color{black}

\vspace{10mm}

\noindent {\bf Acknowledgments.}
The authors are members of the Gruppo Nazionale per l’Analisi Matematica,
la Probabilità e le loro Applicazioni (GNAMPA) of the Istituto Nazionale di Alta Matematica (INdAM). 
A. Passarelli di Napoli has been partially supported through the INdAM-GNAMPA Project 2025 (CUP: E5324001950001) ``Regolarità di soluzioni di equazioni paraboliche a crescita nonstandard degeneri''.
 M. Eleuteri and  A. Passarelli di Napoli have been partially supported through the INdAM-GNAMPA Project 2026 (CUP E53C25002010001) ``Esistenza e regolarità per soluzioni di equazioni ellittiche e paraboliche anisotrope''.

\end{document}